\documentclass[journal]{new-aiaa}
\usepackage[utf8]{inputenc}
\usepackage{textcomp}
\usepackage{multirow}

\usepackage{graphicx}
\usepackage{amsmath}
\usepackage[version=4]{mhchem}
\usepackage{siunitx}
\usepackage{longtable,tabularx}
\usepackage{booktabs}
\makeatletter
\def\blfootnote{\gdef\@thefnmark{}\@footnotetext}
\makeatother

\title{Optimal Replenishment Strategy for Satellite Constellation with Dual Supply Modes}

\author{Jaewoo Kim \footnote{Ph.D. Candidate, Department of Aerospace Engineering, 291 Daehak-Ro.}}
\author{Jaemyung Ahn \footnote{Professor of Aerospace Engineering, 291 Daehak Ro. AIAA Associate Fellow.}}
\affil{Korea Advanced Institute of Science and Technology, Daejeon 34141, Republic of Korea}
\author{Taehyun Sung \footnote{Senior Researcher of Agency for Defense Development, Daejeon 34186, Republic of Korea; thsung@add.re.kr. (Corresponding Author).}}
\affil{Korea Advanced Institute of Science and Technology, Daejeon 34141, Republic of Korea}

\begin{document}

\blfootnote{The contents of this work are conducted as a part of the doctoral thesis of the corresponding author.}

\maketitle

\begin{abstract}
This paper proposes a novel inventory management model for the replenishment strategy of a satellite mega-constellation incorporating dual supply modes: normal and auxiliary. The proposed framework employs an indirect channel for normal supply, wherein spare satellites are initially injected into a parking orbit before transferring to the target orbital plane via propulsion systems and orbital perturbations. Conversely, the auxiliary supply mode utilizes a direct channel, injecting spare satellites immediately into their designated orbital planes. The inventory management model is constructed using parametric replenishment policies: $\left(s,Q\right)$ and $\left(R_1,R_2,Q_1,Q_2\right)$ with a time window. Following this model, two optimization problems are addressed to construct the supply chain for satellite mega-constellation replenishment and to evaluate its performance. These are decision-making contexts from the perspectives of constellation operators and launch service providers. The case study showcases the practical applicability of the proposed model and optimization problems, yielding valuable insights for stakeholders in the spaceborne industry.
\end{abstract}

\section*{Nomenclature}
{\renewcommand\arraystretch{1.0}
\noindent\begin{longtable*}{@{}l @{\quad=\quad} l@{}}
$\lambda_{\text{sat}}$ & failure rate per satellite, in units of satellites per year\\
$\lambda_{\text{plane}}$ & failure rate at in-plane spare, in units of satellites per time unit\\
$\lambda_{\text{parking}}$ & failure rate at parking spare, in units of batches $Q_1$ per time unit\\
$N_t$ & number of discretized time units per year \\
$N_{\text{plane}}$ & number of orbital planes \\
$N_{\text{parking}}$ & number of parking orbits \\
$N_{\text{sats}}$ & number of operational satellites per orbital plane in the constellation \\
$\mu_{\text{primary}}$ & mean waiting time to launch the primary vehicle, time units \\
$t_{\text{primary}}$ & order processing to launch the primary vehicle, time units \\
$\mu_{\text{auxiliary}}$ & mean waiting time to launch the auxiliary vehicle, time units \\
$t_{\text{auxiliary}}$ & order processing time to launch the auxiliary vehicle, time units \\
$R_1$ & reorder point at in-plane spare for indirect channel, satellites\\
$Q_1$ & order quantity at in-plane spare through an indirect channel, satellites \\
$R_2$ & reorder point at in-plane spare for direct channel, satellites \\
$Q_2$ & order quantity at in-plane spare through a direct channel, satellites\\
$Q_{\text{plane}}$ & expected replenishment order quantity at in-plane spare, satellites\\
$k_{R_1,3}$ & reorder point at parking spare, in units of batches $Q_1$\\
$k_{Q_1,3}$ & order quantity at parking spare, in units of batches $Q_1$\\
$p_1$ & probability of one-order cycle \\
$p_2$ & probability of two-order cycle \\
$t_w$ & time window to trigger the second order (direct injection), time units\\
$P_{\text{av}}\left(j^{\text{th}}\right)$ & probability of getting supply from the $j^{\text{th}}$ closest parking orbit \\
$T\left(j^{\text{th}}\right)$ & domain of the Pav($i^{\text{th}}$) \\
$f$ & conditional distribution function of demand (failures) given a time interval\\
$l_1$ & lead time distribution from parking spare to in-plane spare \\
$l_2$ & lead time distribution for direct injection (launch from ground to in-plane spare) \\
$l_3$ & lead time distribution for indirect injection (launch from ground to parking spare) \\
$g$ & inter-order time distribution, from reorder point $R_1$ to reorder point $R_2$ \\
$D_y$  & support of lead time distribution from parking spare to in-plane spare\\
$D_z$  & support of lead time distribution of the auxiliary vehicle\\ 
$D_g$  & support of inter-order time distribution \\
$D_w$  & support of lead time distribution of the primary vehicle\\
$F_{\lambda}$ & cumulative Poisson distribution with parameter $\lambda$ \\
$ES_{\text{plane}}$ & expected shortage at in-plane spare, in units of satellites per cycle\\
$ES_{\text{parking}}$ & expected shortage at parking spare, in units of batches $Q_1$ per cycle\\
$\overline{SL_{\text{plane}}}$ & mean stock level at in-plane spare, in units of satellites\\
$\overline{SL_{\text{parking}}}$ & mean stock level at parking spare, in units of batches $Q_1$ \\
$CS_{\text{plane}}$ & average cycle stock at in-plane spares, in units of satellites times time units\\
$t_{c,\text{plane}}$ & expected replenishment cycle time at in-plane spare, time units \\
$\rho_{\text{plane}}$ & order fill rate at in-plane spare \\
$\rho_{\text{parking}}$ & order fill rate at parking spare \\
$\eta$ & minimum threshold for relative usage of auxiliary launcher \\
$\mathbb{E}$ & expected value of a random variable \\
$\exp\left(\mu\right)$ & exponential distribution with parameter $\mu$
\end{longtable*}}

\section{Introduction}
\lettrine{O}ver the past few decades, numerous satellite constellation missions (including those for Earth observation, communications, and Earth science) have been proposed and successfully operated. Historically, such missions have relied on several costly, high-performance satellites. However, electronics miniaturization, reduced manufacturing and launch costs have enabled developing and deploying new constellations composed of hundreds to thousands of small satellites to meet the increasing demand for global services \cite{wekerle2017status}. The \textit{mega-constellation} initiatives include projects such as OneWeb, which initially proposed a constellation of 6,372 satellites but has since successfully deployed 618 first-generation satellites \cite{pachler2021updated,2023successful}; SpaceX's Starlink program envisioning an ultimate configuration of 42,000 satellites in orbit \cite{iemole2021spacex}, with over 6,000 already launched; and Amazon's Project Kuiper aiming to establish a network comprising 3,236 satellites \cite{pachler2021updated,henry2019amazon}, with successful test of prototype satellites \cite{2024amazon}.

Since a satellite mega-constellation consists of an enormous number of satellites, the failures that the operator must respond to are more frequent than those of a traditional small-sized satellite constellation. In addition to the scale, the preference for reduced reliability to achieve cost-effectiveness makes maintenance more challenging. For example, the design lifetime of a Starlink satellite is 5 years, shorter than the 10 to 20-year lifespans of traditional heavy satellites \cite{iemole2021spacex}. The mega-scale and less reliable satellites of mega-constellations make traditional approaches, such as launching new satellites to replace failed ones on demand or maintaining spare satellites for each orbital plane, limited in their applicability. These new challenges necessitate that operators re-formulate their maintenance strategy.

For the satellite constellation maintenance, the replenishment option that includes the replacement and spare strategies is proposed and studied in the previous work. Formulating the replenishment strategy involves considering where to place, how to interconnect, and how to transport the spares: the \textit{logistics} of spares. The operator can place the spares in orbit (active or standby), in parking, or on the ground \cite{cornara1999satellite}. In these options, the replacement time varies from a day to a year. For example, the replacement time for the in-orbit spare option is approximately 1-2 days since it can immediately replace the failed satellites from the slight difference in altitude relative to the operational orbit. However, the in-parking spare option can replace the failed satellites after several months since it consists of the spare satellites in a lower altitude orbit at the same inclination as the constellation, which requires a transfer maneuver to the operational orbit. Finally, the ground-spare option, always available to replace failed satellites, requires a replacement time of approximately a few months to one year. These spare options can be utilized individually or in combination.

In addition to where to locate the spares, how to interconnect and transport them should also be considered in formulating the logistics. For example, ground spares can be injected into the target orbit by a launch vehicle directly and wait for in-orbit failures as in-orbit spares. However, they can also be initially injected into a parking orbit as parking spares and reach the target orbit indirectly by using an orbital perturbation and propulsion system equipped by the satellites. Decisions on these options form the baseline of the replenishment strategy, followed by the replenishment policy, and thus determine the overall efficiency of the maintenance system.

This research proposes a novel framework for maintaining satellite mega-constellations with dual supply modes. Specifically, this proposed approach uses inventory control methods to assess replenishment strategies in the preliminary phase of satellite mega-constellation projects, which simultaneously uses direct and indirect injection channels. Here, the direct channel injects spare satellite(s) directly into the orbital plane of the constellation via launch service. In contrast, the indirect channel initially places spare satellite(s) into a parking orbit and transfers them into the orbital plane through orbital perturbation and the satellite's propulsion system.

The rest of this paper is organized into five sections. Section II summarizes the literature review and identifies the existing research gaps. Section III provides an appropriate description of the concepts to understand satellite constellations and the fundamentals necessary to express the mathematical model. Section IV develops the mathematical model based on a multi-echelon dual-sourcing inventory management methodology and its validation. Section V describes the relevant optimization problems for satellite mega-constellation maintenance from two different perspectives: the constellation operator and the launch service provider. Section VI covers case studies proposed in Section V and finally, Section VII concludes this paper.

\section{Literature Review}
This section locates our research within the existing literature on the maintenance of satellite constellations. Despite the limited literature on this topic, various approaches have been explored to develop solutions. Notably, simulation-based approaches \cite{cornara1999satellite, lang1998comparison, lansard1998satellite, palmade1998skybridge, ereau1996modeling}, Markovian models \cite{sumter2003optimal, kelley2004minimizing}, and inventory management methods \cite{dishon1966communications, jakob2019optimal} have been applied to replenishment strategy.

Simulation-based approaches were explored by Cornara et al. \cite{cornara1999satellite}, Lang and Adams \cite{lang1998comparison}, Lansard and Palmade \cite{lansard1998satellite}, and Palmade et al. \cite{palmade1998skybridge}. These studies proposed various constellation designs and analyzed replacement strategies, considering options such as on-ground spares, in-parking spares, in-orbit spares, and overpopulation. However, they did not address the selection of mixed strategies, and their criteria for choosing specific strategies were ambiguous.

Ereau and Saleman \cite{ereau1996modeling} adopted a simulation-based approach from an analytic perspective, proposing a constellation deployment and maintenance model using \textit{Petri nets}. They viewed the problem as a non-Markovian distributed system, but their approach faced state explosion issues when incorporating the temporal dimension.

Markovian models were proposed by Sumter \cite{sumter2003optimal} and Kelly and Dessouky \cite{kelley2004minimizing}. Sumter developed an analytic model to find the optimal satellite constellation replacement policy that minimizes maintenance costs using a \textit{Markov Decision Process} (MDP) with a finite state space. The problem-specific costs were linked to satellite purchases, launches, storage, and operational capabilities. Kelly and Dessouky \cite{kelley2004minimizing} used a \textit{Markov chain} to evaluate life cycle costs, employing metaheuristic optimization algorithms like the genetic algorithm and the simulated annealing to find cost-optimal solutions. However, these models have limitations for the large constellations due to exploding dimensions.

Dishon and Weiss \cite{dishon1966communications} and Jakob et al. \cite{jakob2019optimal} explored inventory management methods. Dishon and Weiss modeled the satellite constellation replacement problem using an $\left(s,S\right)$ policy in inventory management theory, analyzing it from a satellite-level perspective. This approach did not account for batch ordering or parking orbit utilization, limiting its effectiveness for mega-constellations.

Jakob et al. \cite{jakob2019optimal} proposed a multi-echelon inventory control for the maintenance strategy of satellite mega-constellations. They successfully integrated parking orbits as part of the supply chain for the replenishment strategy, where the altitude is lower than the operational one, thus the drifting of the parking orbit's right ascension of the ascending node (RAAN) is accounted. They viewed the ground spare as a supplier, parking orbits as shared warehouses, and in-plane spares as retailers. In the research, the optimization formulation is introduced to identify the cost-optimal $\left(s,Q\right)$ policy to minimize the total relevant cost per year for both in-parking and in-plane spares. The total cost comprises the manufacturing, ordering, holding, and maneuvering costs. Sung and Ahn \cite{sung2023optimal} proposed the optimal deployment of satellite mega-constellation problem as expanding the supply chain from single to dual channel by adopting the multi-echelon inventory management model.

Although there are few attempts to design a satellite constellation replenishment policy using inventory management, inventory control, and supply chain management are extensively covered as general topics in the literature. Multi-echelon systems, relevant to satellite spare strategy design, have been explored by Kim \cite{kim1991modeling}, Ganeshan \cite{ganeshan1999managing}, and Schwarz et al. \cite{schwarz1985fill}. Their models featured one warehouse and multiple retailer distribution systems, establishing cost-optimal inventory policies based on expected shortages, order fill rates, and average stock levels. Costantino et al. \cite{costantino2013multi} applied similar models in the aeronautical industry, constructing supply chains dependent on indirect channels.

Dual-sourcing or dual-channel models were developed by Ramasesh \cite{ramasesh1988single}, Ramasesh et al. \cite{ramasesh1991sole}, Chiang and Benton \cite{chiang1994sole}, Balakrishnan \cite{balakrishnan1994order}, Pan et al. \cite{pan1991multiple}, and Schimpel \cite{schimpel2014dual}. These models considered stochastic lead times and constant or normally distributed demand, finding cost-optimal policies through analytic expressions. Especially, Schimpel \cite{schimpel2014dual} expanded the traditional $\left(s,Q\right)$ policy to dual-sourcing models, considering normal and emergency supply modes.

\section{Preliminaries}
This section outlines the fundamentals for developing the mathematical model that this research covers. First, it introduces satellite constellations and orbital mechanics, and the inventory management methodology follows to explain the concepts for formulating the replenishment policy.

\subsection{Satellite Constellation and Orbital Mechanics}
This subsection elaborates on the theory of satellite constellations, orbit perturbations, orbital transfers, and the rocket equation \cite{wertz2001mission, edberg2020design, vallado2022fundamentals}. It summarizes the key concepts required to comprehend this study.

\subsubsection{Satellite Constellation}
A constellation is ``a set of satellites distributed over space to work together to achieve common objectives'' \cite{wertz2001mission}. Notably, constellations exhibit various patterns (configurations) and analysis methodologies.

The Walker Delta pattern is a well-recognized satellite configuration noted for its symmetry. This configuration includes $T$ satellites, distributed evenly among $P$ orbital planes, each containing $S$ satellites. All orbital planes share the same inclination angle $i$. The ascending nodes of these planes are evenly spaced along the equator, separated by intervals of $2\pi/P$. Within each plane, the $S$ satellites are uniformly spaced at intervals of $2\pi/S$. The phase difference in satellites between the adjacent planes is defined by an integer multiple $F$ of $2\pi/T$, ensuring a consistent phase relationship across all planes. This configuration is typically expressed as $i:T/P/F$.

\begin{figure}[hbt!]
    \centering
    \includegraphics[width=0.8\linewidth, page=1]{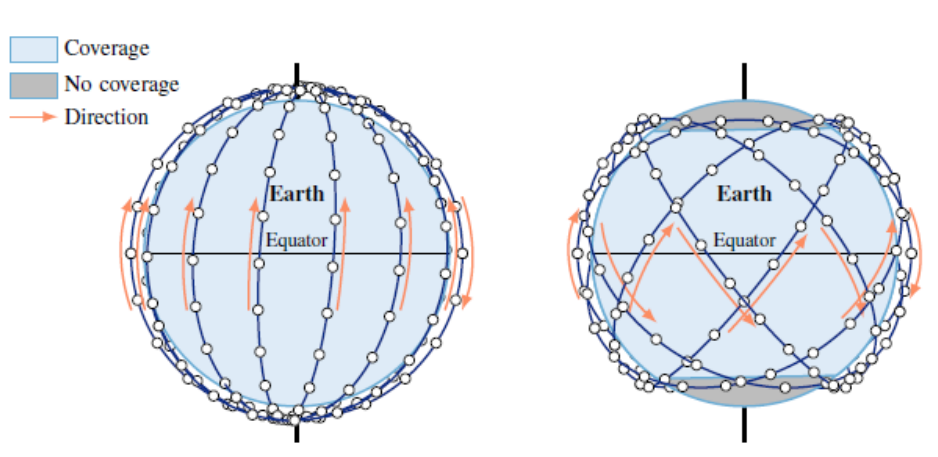}
    \caption{Diagram of Walker Star (left) and Walker Delta (right) geometry \cite{leyva2022ngso}}
    \label{fig: walker patterns}
\end{figure}

The Walker Star pattern, another prominent configuration, is a variation of the Walker Delta pattern. It primarily differs in the spacing of the ascending nodes, set at intervals of $\pi/P$. The Walker pattern ensures that all satellites experience similar orbital perturbations and maintain the same rate of RAAN drift, owing to its symmetricity. In this paper, the satellite constellation's pattern is assumed as Walker pattern, and Figure \ref{fig: walker patterns} illustrates its configuration.

\subsubsection{Orbital Perturbation}
The two primary perturbations influencing satellite motion are the Earth's oblateness and atmospheric drag. This study only considers the Earth's oblateness as the perturbing factor affecting satellite motion. The Earth's oblateness induces secular changes in three orbital elements: the RAAN, the argument of perigee, and the mean anomaly. The primary orbital disturbance relevant to this study is the RAAN drift, where the orbital plane gradually shifts westward/eastward over time at a specific rate. This nodal precession is related to the zonal harmonics coefficient ($J_2$), semi-major axis ($a$), inclination ($i$), and eccentricity of the orbit ($e$) \cite{vallado2022fundamentals}.

\begin{equation}\label{eqn: nodal precession}
    \dot{\Omega} =-\frac{3}{2}  \frac{n  R_E^2}{a^2  (1-e^2)^2}  J_2  \cos{i}
\end{equation}
where $R_E$ is the Earth's equatorial radius. Here, $n$ is the satellite's mean motion, defined as follows, where $\mu_E$ is the Earth's gravitational constant and $a$ is the semi-major axis, respectively.
\begin{equation}\label{eqn: mean motion}
    n=\sqrt{\frac{\mu_E}{a^3}}
\end{equation}

\subsubsection{Orbital Transfer}
Orbital transfer repositions a satellite from one orbit to another. This paper emphasizes fuel-efficient transfers between coplanar circular orbits within an inverse square gravitational field, since most large constellations in the space industry consist of circular orbits.

While impulsive maneuver effectively approximates most maneuvers, this study adopts finite continuous-thrust maneuvers. Low-thrust propulsion systems, such as electric propulsion, have been proposed and applied for continuous-thrust maneuvers of small satellites due to their high efficiency. The required velocity increment for a continuous low-thrust transfer between circular orbits was modeled by Alfano et al. \cite{alfano1994constant}. Bettinger et al. \cite{bettinger2014mathematical} derived the equivalence of velocity increment between Hohmann and continuous low-thrust transfers by comparing them through consecutive Hohmann differences. The required velocity increment ($\Delta v$) and the time-of-flight ($TOF$) for the continuous low-thrust maneuver are as follows:
\begin{equation}\label{eqn: delta v}
    \Delta v = \sqrt{\frac{\mu_E}{r_i}} - \sqrt{\frac{\mu_{E}}{r_f}}
\end{equation}
\begin{equation}\label{eqn: time of flight}
    TOF = \frac{m_{\text{fuel}}}{\dot{m}}
\end{equation}
In Eq. \eqref{eqn: delta v}, $r_i$ and $r_f$ are the radii of the initial and final orbits, respectively, and $\mu_E$ is the Earth's gravitational parameter. In Eq. \eqref{eqn: time of flight}, $m_{\text{fuel}}$ is the required fuel mass from the initial orbit to the final orbit, and $\dot{m}$ is the mass flow rate of the propulsion system.

The following rocket equation relates the velocity increment ($\Delta v$) to the required fuel mass ($m_{\text{fuel}}$):
\begin{equation}\label{eqn: rocket equation}
    m_{\text{fuel}} = m_{\text{dry}}(e^{\Delta v/v_{\text{ex}}}-1)
\end{equation}
Here, $m_{\text{dry}}$ and $v_{\text{ex}}$ are the satellite's dry mass and the exhaust velocity of the propulsion system, respectively.

\subsection{Inventory Management}
Inventory management is a systematic approach to sourcing, storing, and selling inventory, including raw materials (components) and finished goods (products). In business terms, inventory management ensures that the appropriate stock is available at the correct level, at the right time, and at a minimal cost.

This subsection outlines inventory management fundamentals with stochastic lead time and demand. First, we introduce the $\left(s,Q\right)$ policy, one of the most basic parametric replenishment policies. Subsequently, the key parameters necessary for formulating the inventory policy, such as expected shortage (or expected number of backorders) and mean stock level, are discussed. Finally, the $\left(R_1,R_2,Q_1,Q_2\right)$ replenishment policy for dual sourcing, an extension of the $\left(s,Q\right)$ policy, is presented.

\subsubsection{\texorpdfstring{$\left(s,Q\right)$}{} Policy}
The $\left(s,Q\right)$ policy is a continuous review system. This policy orders a fixed quantity $Q$ whenever the inventory position drops to the reorder point $s$. The advantage of the fixed order quantity is its simplicity, unlike other policies such as $\left(R,S\right)$ or $\left(s,S\right)$ \cite{silver2016inventory}. Following the $\left(s,Q\right)$ policy, we can analyze the behavior of the stocks within its cycle, from a reorder point to the next reorder point, which is illustrated in Fig. \ref{fig: (s,Q) policy}.

\begin{figure}[hbt!]
    \centering
    \includegraphics[width=0.8\linewidth, page=2]{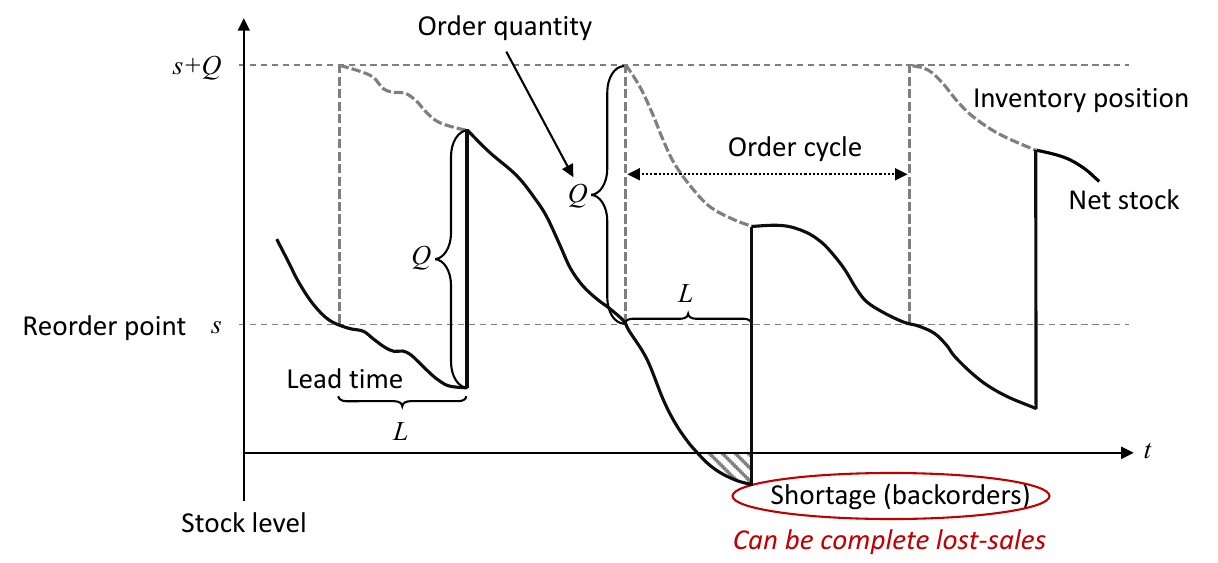}
    \caption{Illustration of the $\left(s,Q\right)$ policy in stock cycle representation}
    \label{fig: (s,Q) policy}
\end{figure}

The following three categories of inventory are defined to describe the system's behavior \cite{silver2016inventory}.
\begin{itemize}
    \item \textit{On-hand stock} ($OH$): This refers to the stock physically present on the shelf, which can never be negative. It determines whether customer demand can be immediately fulfilled.
    \item \textit{Net stock} ($NS$): It is defined as
    \begin{equation}
        \text{Net stock} = (\text{On-hand stock}) - (\text{Backorders})
    \end{equation}
    This quantity can be negative when there are backorders. It is used in various mathematical derivations and is a part of the following crucial definition.
    \item \textit{Inventory position}: It is also known as \textit{available stock}. Inventory position is defined as (Inventory position) = (On-hand stock) + (On-order) - (Backorders) - (Committed). The on-order stock represents the stock that has been requested but not yet received by the stocking point under consideration.
\end{itemize}

Minimizing the relevant inventory cost is the primary objective when determining the inventory policy. The inventory-related costs generally include ordering, holding, and shortage costs. However, estimating the shortage cost is often challenging in practice. Thus, it is preferable to specify a required order fill rate (or service level) to achieve, thereby minimizing the total cost of ordering and holding inventory. The order fill rate ($\rho$) is the fraction of demand that can be immediately satisfied from stock on hand, and the formal definition is given as follows:
\begin{equation}\label{eqn: order fill rate}
    \rho = 1 - \frac{ES}{Q}
\end{equation}
where $ES$ is the expected shortage and $Q$ is the fixed order quantity.

\subsubsection{Expected Shortage}
Evaluating the shortages faced over the replenishment cycle is essential since it determines the service level of the inventory management system. In this study, the service level is represented as the order fill rate, as provided in Eq. \eqref{eqn: order fill rate}. The expected shortage during a lead time is the mean value of the demand exceeding $s$ units from the start of the replenishment cycle until the ordered items arrive.

If the demand $x$ during the replenishment lead time has a probability density function, $f_x$, then the probability that the total demand lies between $x_0$ and $x_0 + dx_0$ is denoted as $f_x(x_0)dx_0$. Thus, the general expression for the expected shortage per replenishment cycle can be presented as follows \cite{silver2016inventory}.
\begin{equation}\label{eqn: expected shortage integral form}
    ES\left(s\right)=\int_{s}^\infty \left(x_0 - s\right)   f_x\left(x_0\right) dx_0
\end{equation}

Specifically, when the demand follows the Poisson process with parameter $\lambda$, the expected shortage during the lead time $\tau$ can be expressed as follows, given that the expectation of demand in the lead time $\tau$ is $\lambda \tau$:
\begin{equation}\label{eqn: expected shortage summation form with Poisson}
    ES\left(s;\lambda,\tau\right)=\sum_{k=s}^\infty \left(k - s\right)  \frac{e^{-\lambda \tau}\left(\lambda \tau\right)^k}{k!}
\end{equation}
Then, Eq. \eqref{eqn: expected shortage summation form with Poisson} can be calculated as Eq. \eqref{eqn: expected shortage expansion with Poisson} \cite{ganeshan1999managing}, where $F_{\lambda \tau}$ represents the cumulative Poisson distribution with parameter $\lambda \tau$.
\begin{equation}\label{eqn: expected shortage expansion with Poisson}
    ES\left(s;\lambda,\tau\right)=\lambda\tau   \left(1 - F_{\lambda \tau}\left(s - 1\right)\right) - s  \left(1 - F_{\lambda \tau}\left(s\right)\right)
\end{equation}
Since this study assumes the satellite's failure follows the Poisson process, Eq. \eqref{eqn: expected shortage expansion with Poisson} is applied to calculate the expected shortage for developing the replenishment policy.

When out of stock, the situation where any demand is back-ordered and filled as soon as adequate replenishment arrives is called \textit{complete back-ordering}. On the contrary, when out of stock, the situation where any demand results in lost sales, because the customer goes elsewhere to satisfy their needs, is called \textit{complete lost sales}. In this study, complete back-ordering is the only concern, and the situation is depicted in Fig. \ref{fig: (s,Q) policy}.

\subsubsection{Mean Stock Level}
The mean stock level is used to approximate the inventory holding cost in this study. Under the assumption that the expected shortage is very small relative to the average on-hand stock, the expected on-hand stock ($\mathbb{E}\left[OH\right]$) can be approximated as the expected net stock ($\mathbb{E}\left[NS\right]$). Thus, the expected on-hand stock right before the ordered stocks arrive ($\mathbb{E}\left[OH_1\right]$) is approximated as follows:
\begin{equation}
    \mathbb{E}\left[OH_1\right] = s - \mathbb{E}\left(x_\tau\right)
\end{equation}
Here, $x_{\tau}$ is the demand during the lead time $\tau$. Since each replenishment is of size $Q$, the expected on-hand stock right after the ordered stocks arrive ($\mathbb{E}\left[\text{OH}_2\right]$) can be expressed as follows:
\begin{equation}
    \mathbb{E}\left[OH_2\right] = \mathbb{E}\left[OH_1\right] + Q = s - \mathbb{E}\left[x_\tau\right] + Q
\end{equation}
Following these relations and assuming that the on-hand stock level drops linearly, the below equation approximates the mean stock level ($\overline{SL}$):
\begin{equation}
    \overline{SL} = \frac{\mathbb{E}\left[OH_1\right] + \mathbb{E}\left[OH_2\right]}{2} + \frac{1}{2} = s - \mathbb{E}\left(x_\tau\right) + \frac{Q}{2} + \frac{1}{2}
\end{equation}
Here, $1/2$ is the continuity correction factor, compensating for the approximation of the linearization of stock level drops to the discrete model of it. If the demand is a Poisson process with parameter $\lambda$, the mean stock level given lead time $\tau$ can be modeled as follows:
\begin{equation}\label{eqn: mean stock level Poisson}
    \overline{SL}\left(s, Q; \lambda,\tau\right)=s-\lambda\tau+\frac{Q}{2}+\frac{1}{2}
\end{equation}

\subsubsection{Stochastic Dual-Sourcing Model: \texorpdfstring{$\left(R_1,R_2,Q_1,Q_2\right)$} Policy}
The stochastic dual-mode replenishment (SDMR) model describes a policy involving two distinct suppliers or supply channels. The policy specifies the reorder properties from each supplier, and the SDMR model used in this research is inspired by the methodology established by Schimpel \cite{schimpel2014dual}.

In the $\left(R_1,R_2,Q_1,Q_2\right)$  policy, a common warehouse where the stock is continuously monitored is concerned. This policy employs two distinct reorder points and order quantities. The primary reorder point, $R_1$, initiates a replenishment order of $Q_1$ units from a supplier, and a secondary reorder point, $R_2$, triggers an order of $Q_2$ units from the other. With the two distinct reorder points, the policy introduces the time window $t_w$, to activate the secondary order. By assessing the inventory level within the $t_w$ period after the first reorder point, the secondary reorder is initiated if the inventory level falls below the point $R_2$. Figure \ref{fig: exemplified cycle of (R1, R2, Q1, Q2) policy} illustrates a representative cycle of the $(R_1, R_2, Q_1, Q_2)$ policy.

\begin{figure}[hbt!]
    \centering
    \includegraphics[width=0.7\linewidth, page=3]{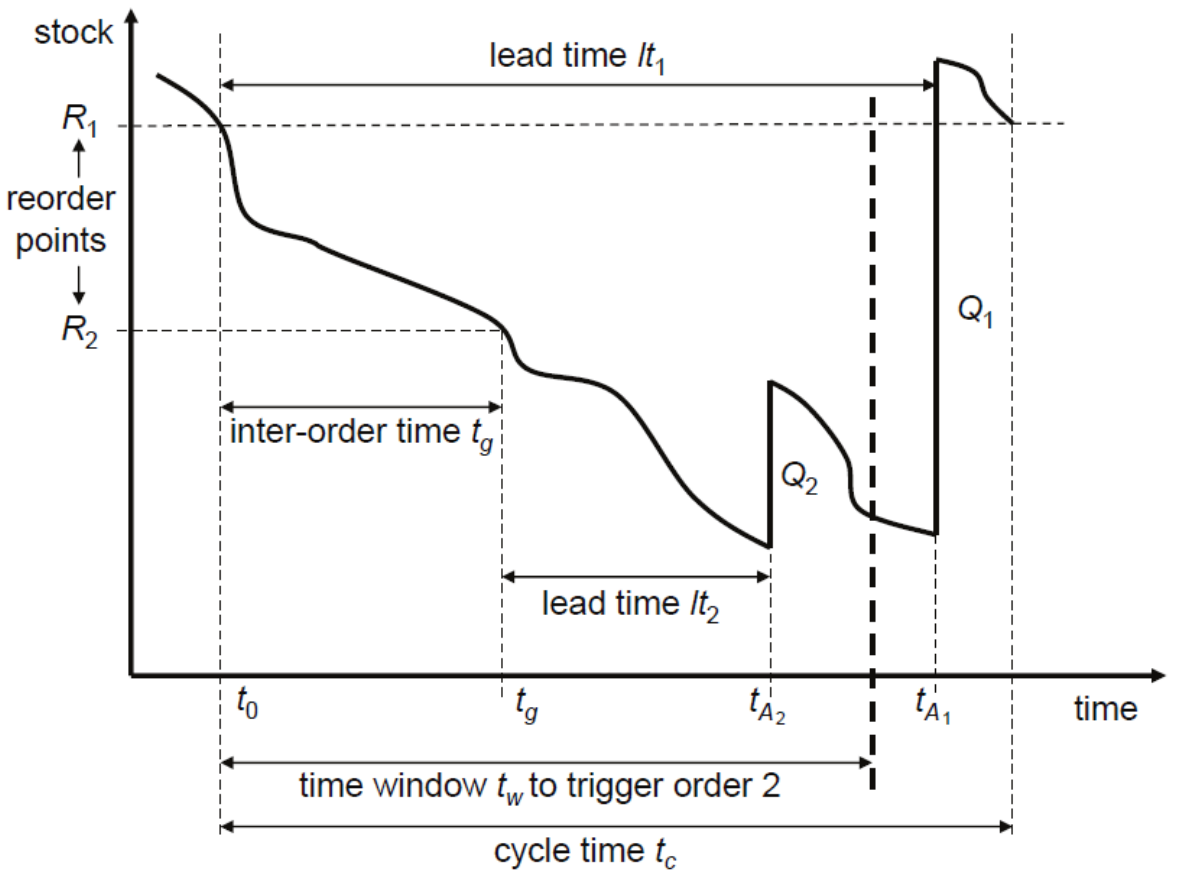}
    \caption{Exemplified cycle of a $\left(R_1, R_2, Q_1, Q_2\right)$ policy with a time window \cite{schimpel2014dual}}
    \label{fig: exemplified cycle of (R1, R2, Q1, Q2) policy}
\end{figure}

If the Possion process models the expenditure of the stock, the demand $x$ given a time interval of the length $\tau$ follows the probability density function $f_{\lambda}(x;\tau)$ with the rate parameter $\lambda$: \begin{equation}\label{eqn: demand conditional lead time}
    f_{\lambda}\left(x;\tau\right)=\frac{\left(\lambda \tau\right)^x}{x!}e^{-\lambda \tau}
\end{equation}
Thus, under the assumption that the failure of the satellite in the orbital plane remains constant, the inter-order time distribution is modeled by the Erlang distribution \cite{kim1991modeling} with the shape parameter $R_1-R_2$ $(R_1>R_2)$ and the satellite's failure rate $\lambda$. The following equation presents the inter-order time ($t$) distribution ($g$):
\begin{equation}\label{eqn: inter-order time}
    g_{\lambda}\left(t;R1,R2\right)=\frac{\lambda^{\left(R_1-R_2\right)}}{\left(R_1-R_2-1\right)!}  t^{\left(R_1-R_2-1\right)}  e^{-\lambda t}
\end{equation}

\begin{figure}[hbt!]
    \centering
    \includegraphics[width=1.0\linewidth, page=4]{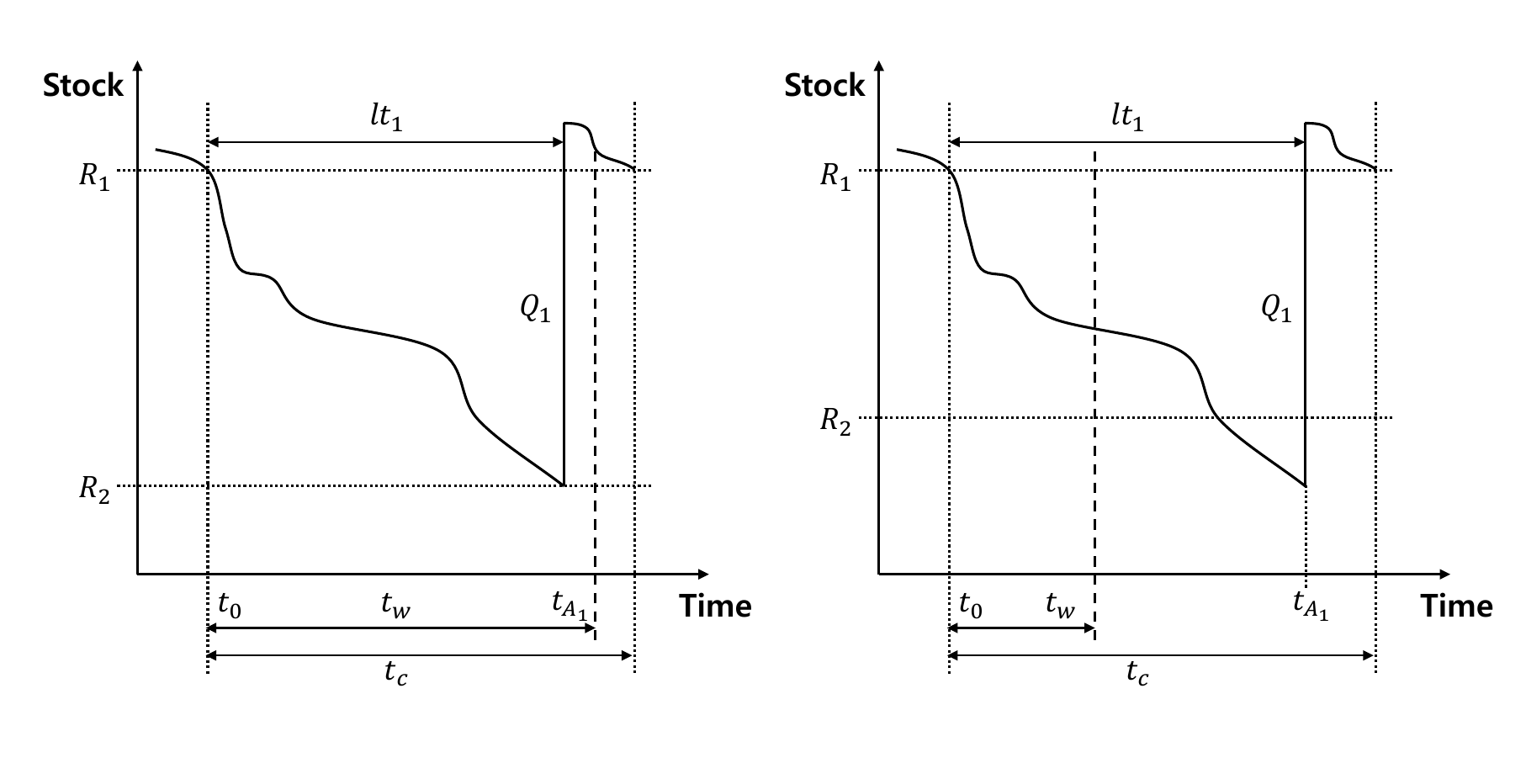}
    \caption{Two possible scenarios of a one-order cycle (order arrival before $t_w$ and after $t_w$)}
    \label{fig: two possible scenarios of a one-order cycle}
\end{figure}

\begin{figure}[hbt!]
    \centering
    \includegraphics[width=1.0\linewidth, page=5]{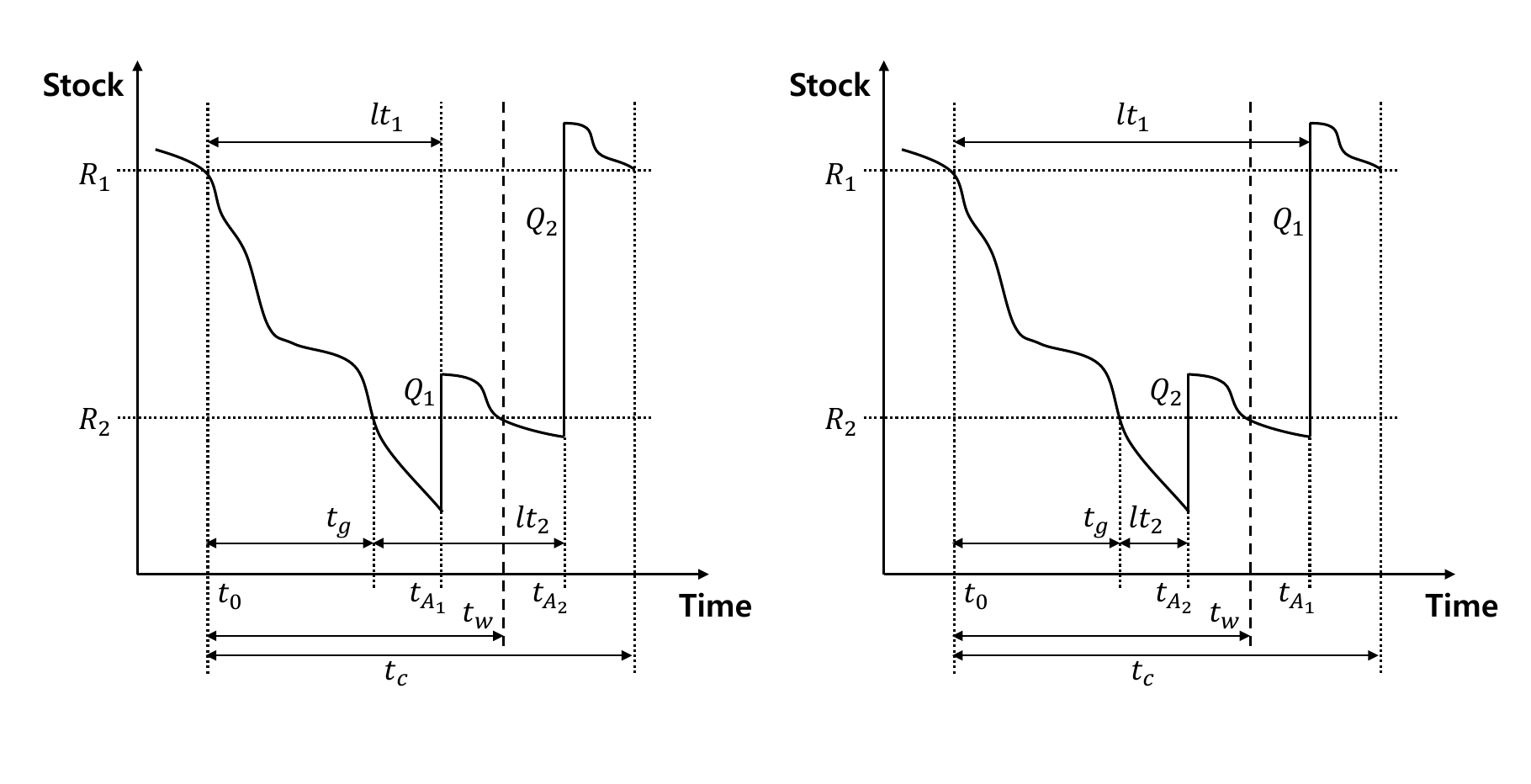}
    \caption{Examples of a two-order cycle ($1^{\text{st}}$ order arrives first, $2^{\text{nd}}$ order arrives first)}
    \label{fig: examples of a two-order cycle}
\end{figure}

Schimpel \cite{schimpel2014dual} proposed an analytic expression for the expected shortage and average stock level in the SDMR model by enumerating all eight possible scenarios, including simultaneous order arrival. There are two cases for the one-order cycle and six cases for the two-order cycle. The one-order cycle is classified by whether the first-order arrival time is longer or shorter than $t_w$. In the one-order cycle cases, the second order should not be triggered during the first-order lead time. For the two-order cycle, the cases are classified by combining the first-order arrival time and the order arrival sequences. The first order can arrive within $t_w$ or not, and the second order can arrive before, after, or at the same time as the first order. The example illustrations of the one-order and two-order cycles are described in Fig. \ref{fig: two possible scenarios of a one-order cycle} and Fig. \ref{fig: examples of a two-order cycle}.

Following this classification, the previous study presented the eight possible cases with their relevant formulae for the dual-sourcing policy. Note that in the continuous time domain, the simultaneous arrival probability of both orders is conceptually zero but is not in the discrete case. Therefore, this paper divides all possible eight scenarios into six cases by including the simultaneous order arrival case into its relevant case. The rearranged cases are summarized below. 

\begin{itemize}
    \item Case 1: One-order cycle, the first order arrives before the time window $t_w$
    \item Case 2: One-order cycle, the first order arrives after the time window $t_w$
    \item Case 3: Two-order cycle, the first order arrives first, before the time window $t_w$
    \item Case 4: Two-order cycle, the first order arrives last, before the time window $t_w$, including simultaneous order arrival
    \item Case 5: Two-order cycle, the first order arrives first, after the time window $t_w$
    \item Case 6: Two-order cycle, the first order arrives last, after the time window $t_w$, including simultaneous order arrival
\end{itemize}

These cases form a partition to cover all the possible events, classified based on the order cycle and arrival sequences. The expected shortage of the SDMR model in this research adheres to the formulas established in the previous study \cite{schimpel2014dual}. Furthermore, the mean stock level is appropriately adjusted to suit the specific context of the proposed problem. 

\section{Model Construction}
This section outlines the inventory management model for the maintenance strategy of satellite mega-constellations with a mixed injection option. Initially, an overview and the assumptions necessary for model construction are presented. Subsequently, the mathematical model for each spare stock is formulated, encompassing the demand (satellite's failure), lead time distribution, expected shortage, and mean stock level throughout a replenishment cycle. The relevant cost model is defined following these inventory models, and finally, the validation of the proposed model concludes this section by assessing its results with those obtained from the Monte Carlo simulation.

\subsection{Overview of the Model and Assumptions}
\begin{figure}[hbt!]
    \centering
    \includegraphics[width=0.6\linewidth, page=6]{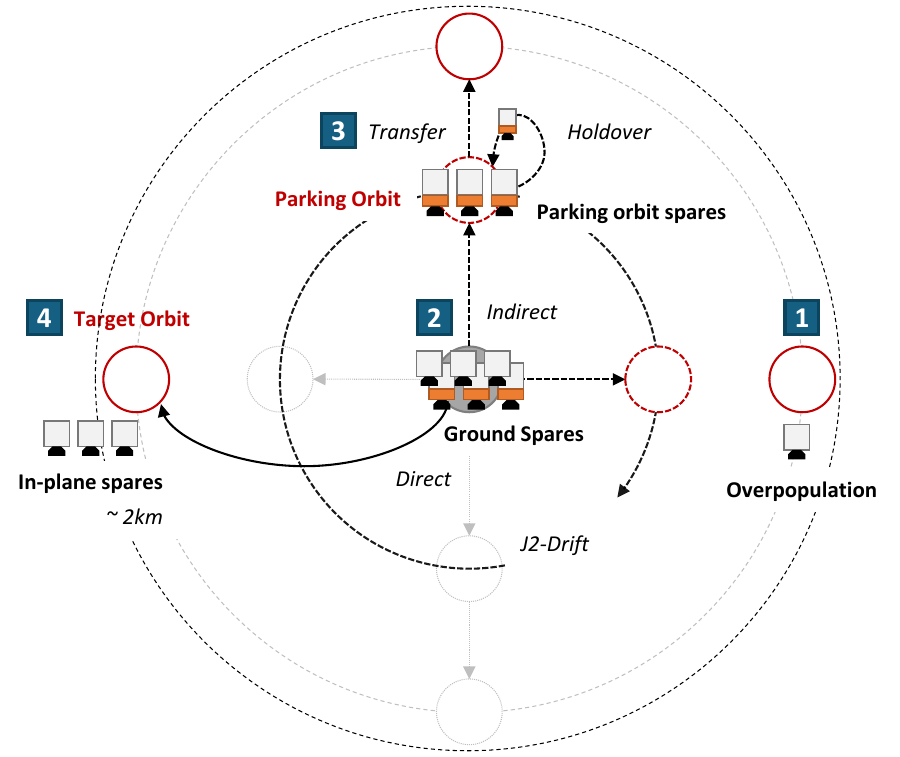}
    \caption{Candidates for replenishment options}
    \label{fig: proposed replenishment options}
\end{figure}

In this paper, the constellation maintenance strategies are developed based on the satellite replenishment option, as presented in the previous research \cite{cornara1999satellite}. Figure \ref{fig: proposed replenishment options} visualizes these options, which are overpopulation, launch-on-demand, in-plane spares, and parking spares.

Among the options, Jakob et al. \cite{jakob2019optimal} proposed a mixed strategy incorporating levels of spares, which serves as the baseline for this study. The proposed strategy can enhance flexibility in designing the spare options for satellite constellations.

The model constructs the first level of spares as the constellation's in-plane spares. The overpopulation and the in-plane spare strategy are not differentiated, assuming failed satellites are replaced directly by either. When a satellites satellite fails, and if an in-plane spare is available, the spare is used to avoid outages with little to no time delay.

The second level of spares is the parking spares, which are located in lower altitude orbits with the same inclination as the orbital planes. If the in-plane spare level drops to a reorder point, parking spares are used to resupply it. Additionally, parking spaces can replenish any orbital plane within the constellation, owing to their relative RAAN drift, thus enhancing the flexibility of the supply chain.

The final level of spares is ground spares, which are assumed to always be available to replenish the orbits. When parking spares reach a reorder level, ground spares are ordered to schedule a launch and replacement. This multi-echelon inventory systems accounts for the stochastic lead times and demand across all spare locations.

In this paper, the sourcing strategy is expanded from single to dual. The baseline model relies on parking orbits for spare satellite injection, which entails orbit transfer maneuvers--a single-sourcing strategy utilizing an indirect channel only. By adopting a direct injection option, ground spares can directly resupply in-plane spares, creating a dual-sourcing strategy where spares are injected in both direct and indirect channels. This proposed mixed injection option can offer additional flexibility and cost advantages, but its complexity complicates to design of its operational criteria.

\begin{figure}[hbt!]
\centering
\includegraphics[width=1.0\linewidth, page=7]{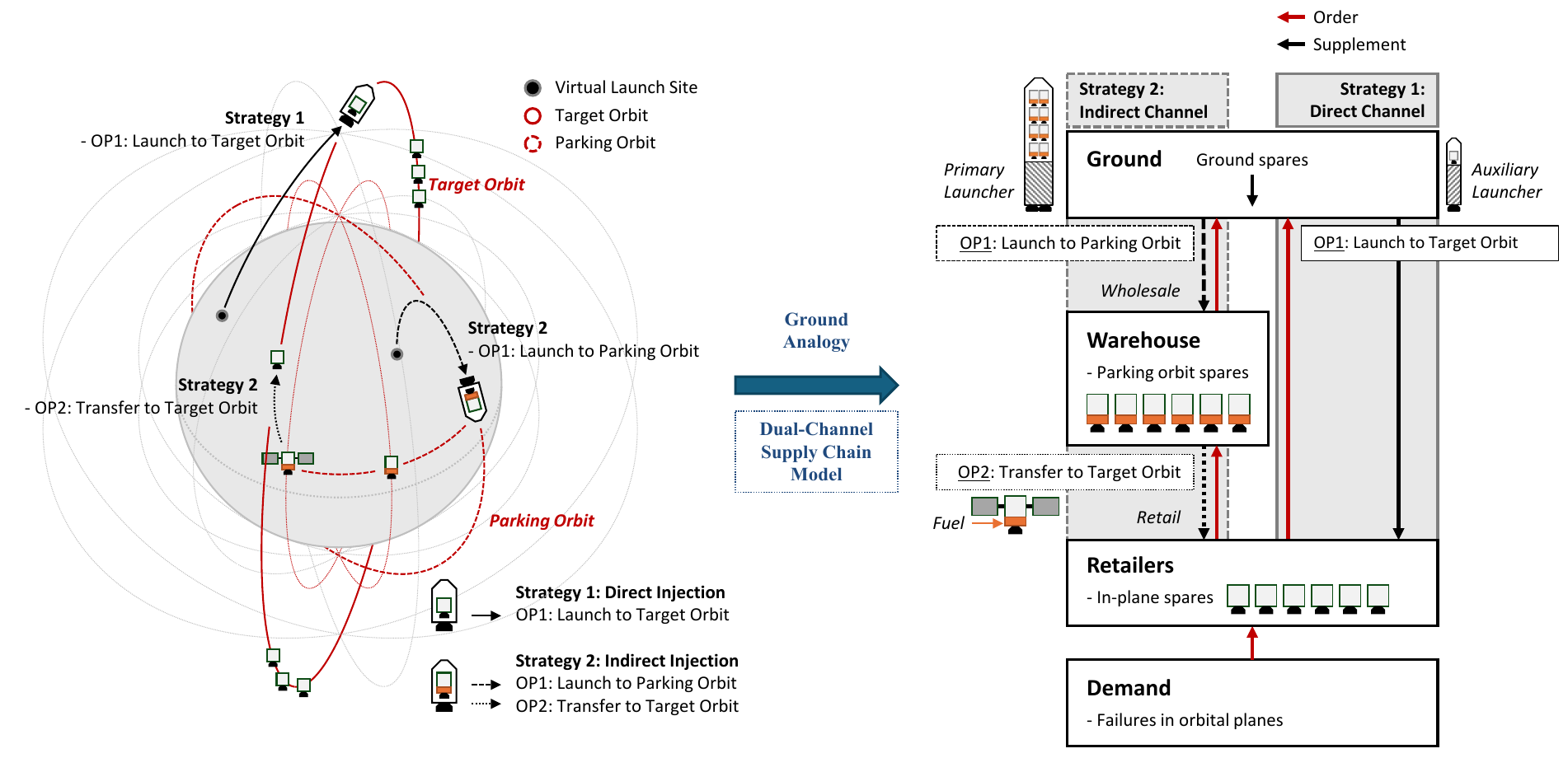}
\caption{Proposed supply chain for satellite mega-constellation replenishment strategy}
\label{fig: proposed supply chain for satellite mega-constellation replenishment strategy}
\end{figure}

\begin{figure}[hbt!]
\centering
\includegraphics[width=0.8\linewidth, page=8]{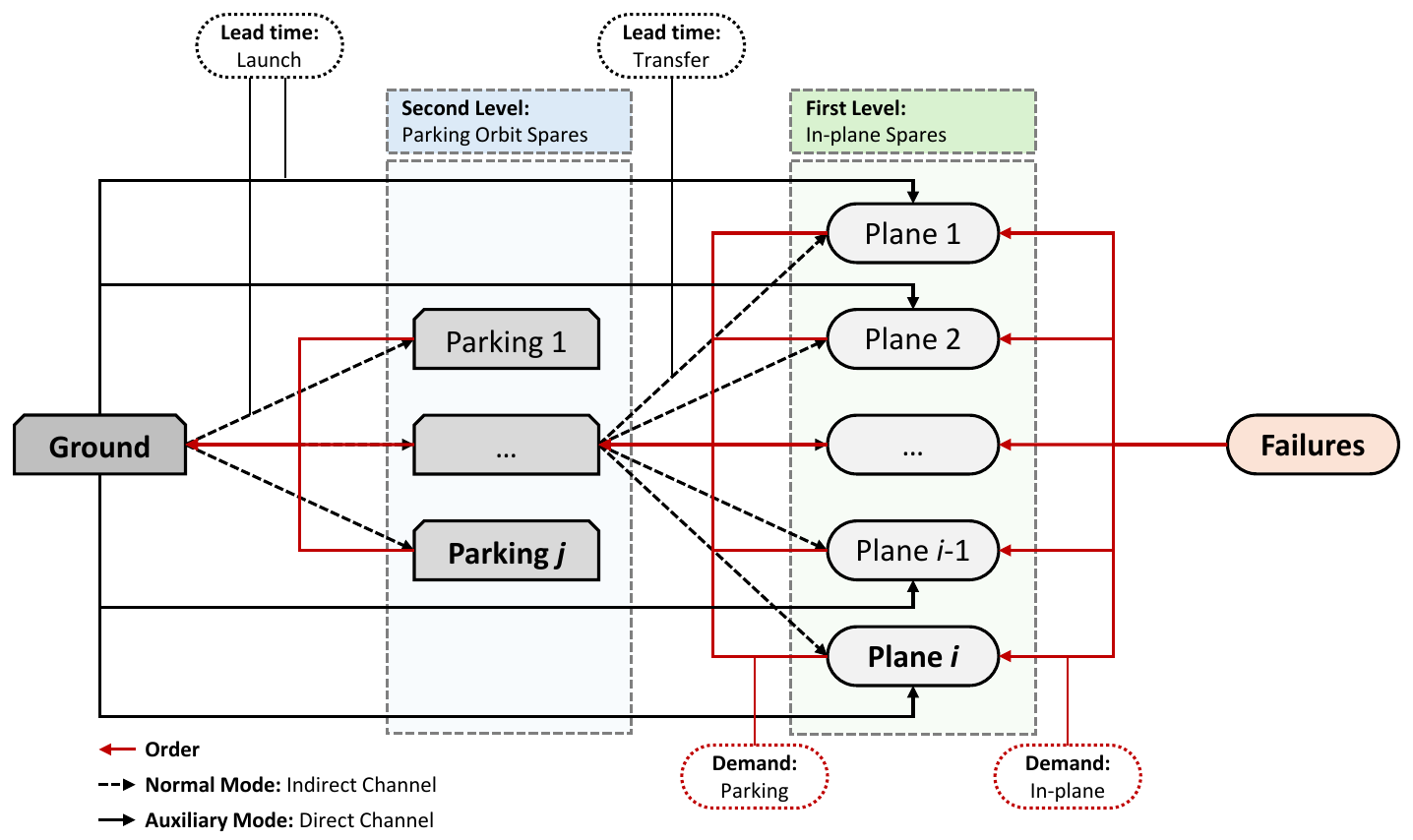}
\caption{Dual-sourcing multi-echelon inventory control model}
\label{fig: dual-sourcing multi-echelon inventory control model}
\end{figure}

Therefore, the parametric $\left(R_1, R_2, Q_1, Q_2\right)$ policy \cite{schimpel2014dual} is employed as a benchmark for the model construction. The proposed model includes normal and auxiliary supply modes, with subscripts 1 and 2 representing each mode, respectively. The indirect channel (normal supply mode) and the direct channel (auxiliary supply mode) are assumed under $R_1 > R_2$ conditions. Primarily, spare satellites are inserted via indirect injection by a primary (heavy) launch vehicle, while any spares not covered by normal operations are inserted by auxiliary (responsive) launchers as needed. Figures \ref{fig: proposed supply chain for satellite mega-constellation replenishment strategy} and \ref{fig: dual-sourcing multi-echelon inventory control model} visually represent the proposed supply chain and the dual-sourcing multi-echelon inventory control model, respectively.

The following list summarizes the general assumptions for the proposed model. The assumptions are discussed as the model is introduced in more detail.

\begin{enumerate}
    \item Spare satellites in the first echelon (in-plane spares) are considered immediate to replace a failed satellite unit. This postulate is true in the case of an overpopulated strategy. The replacement process can take up to 2 days, but this delay is not considered in the model \cite{jakob2019optimal}.
    \item The constellation’s in-plane spares get supplies from the closest available parking orbit’s stock. To allow flexibility in the spare replacement flow, the model allows any parking orbits to potentially resupply any orbital plane’s spares. When a constellation orbit’s in-plane stock level reaches the reorder point, an order is placed to all parking orbits jointly, and the spare batch is supplied from the closet parking orbit with availability at the time of the order \cite{jakob2019optimal}.
    \item The constellation is limited to the Walker pattern \cite{jakob2019optimal}.
    \item Supply from the ground can be delivered only to a unique orbit. Using a single rocket launch to supply different orbital planes can be inefficient \cite{jakob2019optimal}.
    \item An order can only be processed when no previous order is in transit between different spare levels, ensuring only one cycle order is treated at any time. However, multiple launches from the ground simultaneously are permissible.
    \item The expected shortage is very small relative to the average on-hand stocks. Thus, the shortage is negligible for calculating the mean stock levels.
    \item The reorder points $R_1$ and $R_2$ in the dual-sourcing policy can be negative, but the relationship $R_1 > R_2$ must always be maintained. The order quantities $Q_1$ and $Q_2$ are positive integers.
    \item The indirect and direct channels correspond to normal and auxiliary supply modes, respectively. A primary launcher handles the indirect injection, while an auxiliary launcher handles the direct injection. Indices are given to the relevant transportation strategies. For the indirect channel, the orbital transfer from the parking orbits to the orbital planes is denoted as `1,' and the launch to the parking orbits is denoted as `3.' In the case of the direct channel, the launch to the orbital planes is denoted as `2.'
    \item The second order (identifier 2, auxiliary supply mode) is not permitted after the time window $t_w$. Consequently, the auxiliary supply mode can only be triggered if the in-plane spares reach the second reorder point $R_2$ before $t_w$ elapses following the first order.
    \item As the spares have to be transferred by batches both from the Earth’s ground to the parking orbits and from the parking orbits to the constellation orbits, the order quantity and re-order point at the parking orbits are assumed to be multiples of the batch size $Q_1$ of in-plane spares \cite{jakob2019optimal}:
    \begin{equation}
        Q_{\text{parking}} = k_{Q_1,3}  Q_1
    \end{equation}
    \begin{equation}
        R_{\text{parking}} = k_{R_1,3}  Q_1
    \end{equation}
\end{enumerate}

\begin{figure}[hbt!]
\centering
\includegraphics[width=1.0\linewidth, page=9]{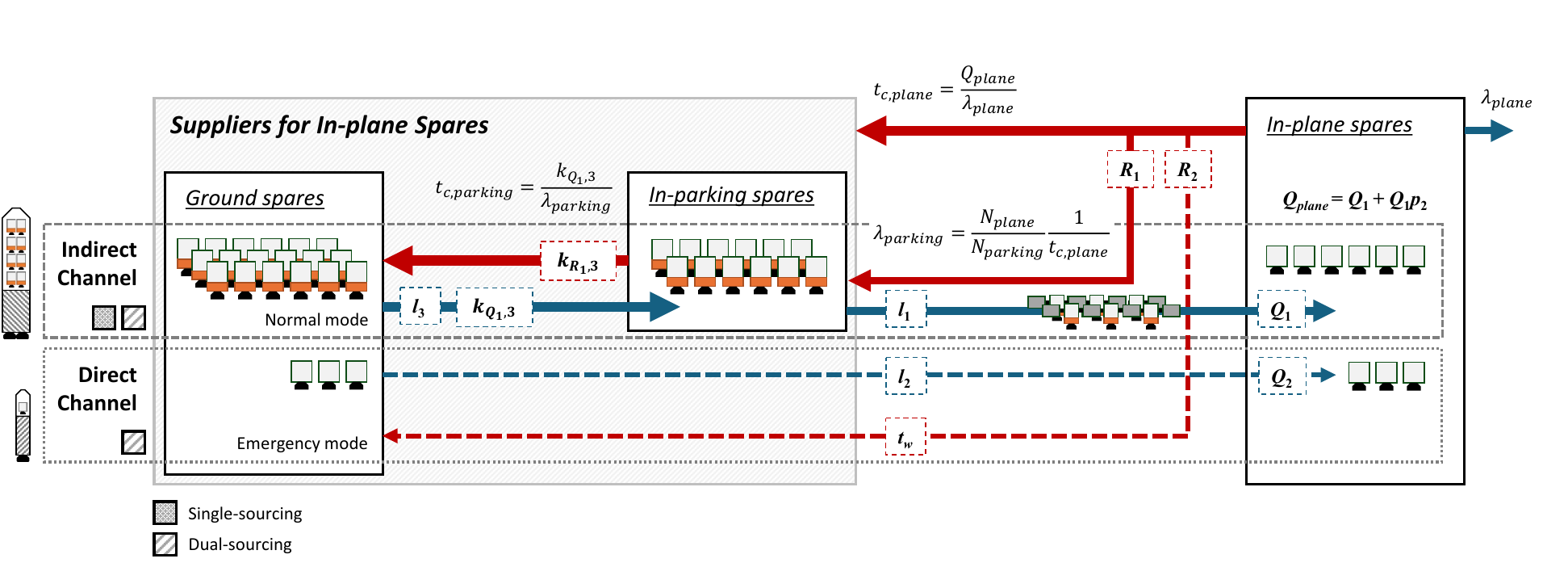}
\caption{Summary of the satellite mega-constellation replenishment model}
\label{fig: Summary of the satellite mega-constellation replenishment model}
\end{figure}

Figure \ref{fig: Summary of the satellite mega-constellation replenishment model} illustrates the graphical representation of the proposed satellite mega-constellation replenishment model. The notations are explained with the model description in the following sections.

\subsection{Inventory Model of In-Plane Spares}
This subsection introduces the inventory model for in-plane spares. The model consists of the dual-sourcing strategy employing an $\left(R_1, R_2, Q_1, Q_2)\right)$ policy with a time window $t_w$ for triggering the second order.

\subsubsection{Demand Model}
The satellite failure is modeled using a Poisson distribution with the parameter representing the satellite's failure rate. The failure rate indicates the number of failures per unit of time, so the satellite's failure rate per constellation orbital plane is expressed as follows\cite{jakob2019optimal}:
\begin{equation}\label{eqn: failure rate plane}
    \lambda_{\text{plane}} = \frac{N_{\text{sats}}  \lambda_{\text{sat}}}{N_t}
\end{equation}
where $N_{\text{sats}}$ is the number of operational satellites per orbital plane, and $\lambda_{\text{sat}}$ is the failure rate of a satellite per year. $N_t$ represents the number of discretized time units per year. For instance, if the base time unit is a day, $N_t$ would be 365, and if the time unit is a week, $N_t$ would be 52. 

\subsubsection{Lead Time Distribution}
Two sourcing channels are available for in-plane spares: indirect and direct, leading to different lead time distributions. First, the lead time distribution for the indirect channel follows the model presented in a previous study \cite{jakob2019optimal}. The multi-echelon spare inventory model leverages the Earth's gravitational field, allowing a parking orbit to visit all nominal constellation planes and supply spares for failures. The transfer maneuver from the parking to nominal orbit is assumed as the in-plane maneuver which performs after the RAAN alignment between them. Consequently, the lead time for transferring batches of satellites through the indirect channel consists of the drift time to align the orbital planes and the time-of-flight of the transfer maneuver.
The probability distribution for transfer time, representing the lead time from parking orbits to the orbital planes, is defined analytically. Assuming the parking constellation follows a Walker pattern, the parking orbits are equally distributed in the equatorial plane, dividing the possible RAAN differences for the drift into $N_{\text{parking}}$ intervals as follows:
\begin{equation}\label{eqn: delta omega}
    \Delta \Omega_j = \left[\left(j-1\right) \frac{2\pi}{N_{\text{parking}}},j  \frac{2\pi}{N_{\text{parking}}}\right),\ \forall j\in \left\{1,2,\dots,N_{\text{parking}}\right\}
\end{equation}
where $\Delta \Omega_j$ is the RAAN difference for $j^{\text{th}}$ closest parking orbit relative to the orbital plane. Since the parking orbits are equally distributed, circular, and share a common radius and inclination, all parking orbits exhibit the same drift rate, as described by Eq. \eqref{eqn: nodal precession}. Therefore, we assume that the alignment time can be expressed as uniformly distributed in each possible interval. The domain of supplement time from $j^{\text{th}}$ closest parking orbit $T_j$is expressed as follows:

\begin{equation}\label{eqn: domain of supplement time}
    T_j=\left[\left(j-1\right) \frac{2\pi}{N_{\text{parking}}}  \frac{1}{\left|\dot{\Omega}_{\text{relative}}\right|}+t_{\text{trans}},j  \frac{2\pi}{N_{\text{parking}}}  \frac{1}{\left|\dot{\Omega}_{\text{relative}}\right|}+t_{\text{trans}}\right),\ \forall j\in \left\{1,2,\dots,N_{\text{parking}}\right\}
\end{equation}
where $t_{\text{trans}}$ is the in-plane transfer maneuver time, and $\left|\Omega\right|$ is the RAAN drift rate of the parking orbit relative to the constellation's orbital plane. 

The concept of the availability of the parking orbit is presented to express the lead time distribution for the transfer. The probability of the parking orbit availability is derived using a binomial-like distribution based on the parking orbit's order fill rate, $\rho_{\text{parking}}$. The probability of getting supply from $j^{\text{th}}$ closest parking orbit $P_{\text{av}}\left(j\right)$ is obtained by:
\begin{equation}
    P_{\text{av}}\left(j\right) = \left(1-\rho_{\text{parking}}\right)^{j-1}  \rho_{\text{parking}},\ \forall j\in\left\{1,2,\dots,N_{\text{parking}}\right\}
\end{equation}

As a result, the lead time distribution of the spare supplement from the parking orbits is expressed as follows:
\begin{equation}\label{eqn: lead time distribution 1}
    l_1\left(y\right)=P_{av,n}\left(j\right)  \frac{N_{\text{parking}}  \left|\dot{\Omega}_{\text{relative}}\right|}{2\pi},\ y\in T_j
\end{equation}
The right-hand side consists of two components. The first component $P_{\text{av},n}\left(j\right)$ is the normalized probability of availability of the $j^{\text{th}}$ parking orbit ($=P_{\text{av}}\left(j\right)/\sum_{k=1}^{N_{\text{parking}}}P_{\text{av}}\left(k\right)$), and the second component is the inverse of interval length of $T_j$. To derive the lead time distribution, it is important to note that the probability of all parking orbits being out-of-stock at the time of resupply is very small thus, can be neglected. Consequently, it is assumed that at least one parking orbit is always available, although it may not necessarily be the closest one. This assumption holds only under the condition of a high order fill rate for each parking orbit.

On the other hand, the direct channel uses auxiliary launchers to resupply the in-plane spares from the ground. The lead time distribution for the direct channel, denoted as $l_2$, includes the launch ordering process and the waiting time for the next launch window, as proposed in previous work \cite{jakob2019optimal}. This study assumes that ground spares are always available. Thus, the availability of ground spares is not considered in modeling the lead time. The order processing time is treated as constant \cite{ganeshan1999managing}, while the waiting time for the next launch window follows an exponential distribution \cite{jakob2019optimal}. Therefore, the lead time for direct injection ($z$) is represented as follows:
\begin{equation}\label{eqn: lead time plane 2}
    z\sim \exp\left(\mu_{\text{auxiliary}}\right)+t_{\text{auxiliary}}
\end{equation}
where $\mu_{\text{auxiliary}}$ is the mean launch interval or waiting time for the next launch of the auxiliary vehicle, and $t_{\text{auxiliary}}$ is the launch order processing time.

\subsubsection{Expected Shortage}
The expected shortage is calculated based on the partitioned sample space, as introduced in Section III.B.4. The values of all the cases are summed to obtain the expected shortages for the entire sample space. The detailed derivation of the original formulation can be found in the previous work \cite{schimpel2014dual}. However, this paper will only introduce the problem-specific formulation.

Consider two order modes with order quantities $Q_1$ and $Q_2$, and their respective lead times $l_1$ and $l_2$. The demand (failures) $x$ over any time interval of length $\tau$ follows the probability density function $f_{\lambda_{\text{plane}}}\left(x;\tau\right)$. Given the two reorder points $R_1$ and $R_2$ for the two order modes, the inter-order time $t$ is a random variable with its associated demand distribution $f_{\lambda_{\text{plane}}}\left(x;\tau\right)$. The model introduces the time window $t_w$ to support deciding whether to trigger the second order. The second order will be triggered if the stock level falls below $R_2$ within the time window $t_w$.

The indices in the subscript imply the cases, and the ordered indices in the superscript define the order arrival sequences. Specifically, $ES^{1}$ denotes a one-order cycle, while $ES^{1,2}$ and $ES^{2,1}$ denote two-order cycles where the first order arrives before the second order, and the second order arrives before the first order, respectively. Additionally, the variables $y$, $t$, and $z$ represent the first-order lead time, inter-order time, and second-order lead time, respectively.

The expected shortage for all possible cases can be expressed with the defined notations. Note that $g\left(t\right)$ substitutes $g_{\lambda_{\text{plane}}}\left(t;R_1,R_2\right)$, and the support for functions $l_1$, $g$, and $l_2$ are denoted as $D_{y}$, $D_{t}$, and $D_{z}$, respectively. Therefore, the following expressions hold for $y \in D_{y}$, $t \in D_{t}$, and $z \in D_{z}$.

\begin{itemize}
    \item Case 1
    \begin{equation}\label{eqn: plane expected shortage 1}
        ES_{1}=\int_{y\le t_w} ES^1\left(y;R_1,R_2,Q_1,Q_2,t_w\right)   l_1\left(y\right) dy
    \end{equation}
    \item Case 2
    \begin{equation}\label{eqn: plane expected shortage 2}
        ES_{2}=\int_{y> t_w} ES^1\left(y;R_1,R_2,Q_1,Q_2,t_w\right)  l_1 \left(y\right)dy
    \end{equation}
    \item Case 3
    \begin{equation}\label{eqn: plane expected shortage 3}
        ES_{3}=\int_{y\le t_w}\int_{t<y}\int_{z>y-t} ES^{1,2}\left(y,t,z;R_1,R_2,Q_1,Q_2,t_w\right)  l_1\left(y\right)  g\left(t\right)  l_2\left(z\right)dydtdz
    \end{equation}
    \item Case 4
    \begin{equation}\label{eqn: plane expected shortage 4}
        ES_{4}=\int_{y\le t_w}\int_{t<y}\int_{z\le y-t} ES^{2,1}\left(y,t,z;R_1,R_2,Q_1,Q_2,t_w\right)  l_1\left(y\right)  g\left(t\right)  l_2\left(z\right)dydtdz
    \end{equation}
    \item Case 5 
    \begin{equation}\label{eqn: plane expected shortage 5}
        ES_{5}=\int_{y> t_w}\int_{t<y}\int_{z>y-t} ES^{1,2}\left(y,t,z;R_1,R_2,Q_1,Q_2,t_w\right)  l_1\left(y\right)  g\left(t\right)  l_2\left(z\right)dydtdz
    \end{equation}
    \item Case 6 
    \begin{equation}\label{eqn: plane expected shortage 6}
        ES_{6}=\int_{y> t_w}\int_{t<y}\int_{z\le y-t} ES^{2,1}\left(y,t,z;R_1,R_2,Q_1,Q_2,t_w\right)  l_1\left(y\right)  g\left(t\right)  l_2\left(z\right)dydtdz
    \end{equation}
\end{itemize}

The detailed model formulations for the inner integral variables are presented in Appendix A. The expected shortage for the in-plane stocks $ES_{\text{plane}}$ can be obtained by summing the expected shortages of all the partitions:
\begin{equation}\label{eqn: expected shortage plane}
    ES_{\text{plane}}=\sum_{k=1}^6 ES_{k}
\end{equation}

\subsubsection{Mean Stock Level}
The calculation of the mean stock level is similar to that of the expected shortage. First, partition the sample space into the six cases introduced in Section III.B.4. Next, calculate the average cycle stock for each case, which is conceptually the integration of the stock level over the cycle time. Then, sum the average cycle stocks of all the cases to obtain $CS_{\text{plane}}$. Finally, divide $CS_{\text{plane}}$ by the expected duration of a replenishment cycle $t_{c,\text{plane}}$ to obtain the mean stock level $\overline{SL_{\text{plane}}}$.

The average cycle stocks for all possible cases are expressed below, following the same subscript and superscript conventions used in deriving expected shortage. Note that the following expressions hold for $y \in D_{y}$, $t \in D_{t}$, and $z \in D_{z}$, which are the supports for the functions $l_1$, $g$, and $l_2$, respectively. Also, $g\left(t\right)$ substitutes $g_{\lambda_{\text{plane}}}\left(t;R_1,R_2\right)$.
\begin{itemize}
    \item Case 1
    \begin{equation}\label{eqn: average cycle stock}
        CS_{1}=\int_{y\le t_w} CS^1\left(y;R_1,R_2,Q_1,Q_2,t_w\right)  l_1\left(y\right) dy
    \end{equation}
    \item Case 2
    \begin{equation}
        CS_{2}=\int_{y> t_w} CS^1\left(y;R_1,R_2,Q_1,Q_2,t_w\right)  l_1\left(y\right) dy
    \end{equation}
    \item Case 3
    \begin{equation}
        CS_{3}=\int_{y\le t_w}\int_{t<y}\int_{z>y-t} CS^{1,2}\left(y,t,z;R_1,R_2,Q_1,Q_2,t_w\right)  l_1\left(y\right)  g\left(t\right)  l_2\left(z\right)dydtdz
    \end{equation}
    \item Case 4
    \begin{equation}
        CS_{4}=\int_{y\le t_w}\int_{t<y}\int_{z\le y-t} CS^{2,1}\left(y,t,z;R_1,R_2,Q_1,Q_2,t_w\right)  l_1\left(y\right)  g\left(t\right)  l_2\left(z\right)dydtdz
    \end{equation}
    \item Case 5 
    \begin{equation}
        CS_{5}=\int_{y> t_w}\int_{t<y}\int_{z>y-t} CS^{1,2}\left(y,t,z;R_1,R_2,Q_1,Q_2,t_w\right)  l_1\left(y\right)  g\left(t\right)  l_2\left(z\right)dydtdz
    \end{equation}
    \item Case 6 
    \begin{equation}
        CS_{6}=\int_{y> t_w}\int_{t<y}\int_{z\le y-t} CS^{2,1}\left(y,t,z;R_1,R_2,Q_1,Q_2,t_w\right)  l_1\left(y\right)  g\left(t\right)  l_2\left(z\right)dydtdz
    \end{equation}
\end{itemize}
The detailed model formulation for the inner integral variables is provided in Appendix B. The average cycle stock for the in-plane spares $CS_{\text{plane}}$ is calculated by summing the average cycle stocks of all the cases:
\begin{equation}
    CS_{\text{plane}}=\sum_{k=1}^6 CS_{k}
\end{equation}

Following this, the mean stock level can be obtained by dividing the average cycle stock $CS_{\text{plane}}$ by the expected cycle time $t_{c,\text{plane}}$:
\begin{equation}
    \overline{SL_{\text{plane}}}=\frac{CS_{\text{plane}}}{t_{c,\text{plane}}}+\frac{1}{2}
\end{equation}
Note that the continuity correction factor $1/2$ is added to adjust the linearization of the stock expenditure. The replenishment cycle $t_{c,\text{plane}}$ is equivalent to the time until the delivered order quantities are depleted by demand (failures). Thus, the expected cycle time is established below with the given failure rate $\lambda_{\text{plane}}$:
\begin{equation}\label{eqn: expected cycle time plane}
    t_{c,\text{plane}} = \frac{Q_{\text{plane}}}{\lambda_{\text{plane}}}
\end{equation}
The $Q_{\text{plane}}$ is the expected reorder quantity for in-plane spares. Since the model adopts a dual-sourcing strategy, the expected reorder quantity can vary depending on the probability of one-order and two-order cycles. Suppose the probability of the one-order cycle is $p_1$, and the probability of the two-order cycle is $p_2$, the expected reorder quantity for in-plane spares is given below:
\begin{equation}\label{eqn: expected reorder quantity plane}
    Q_{\text{plane}} = Q_1   p_1 + \left(Q_1 + Q_2\right)   p_2 = Q_1 + Q_2   p_2
\end{equation}
The expression for the probability of order cycles is summarized in Appendix C. The mean stock level formulation for the in-plane spare is obtained as:
\begin{equation}\label{eqn: mean stock level plane}
    \overline{SL_{\text{plane}}}=CS_{\text{plane}}   \frac{\lambda_{\text{plane}}}{Q_1 + Q_2 p_2} + \frac{1}{2}
\end{equation}

\subsection{Inventory Model of Parking Spares}
This subsection introduces the inventory model for parking spares, employing the $\left(s,Q\right)$ replenishment policy. Note that throughout this study, the $\left(s,Q\right)$ policy is denoted as $\left(R,Q\right)$ for consistent notations.

\subsubsection{Demand Model}
In the ordering process from the operational plane, an order is placed for every $Q_{\text{plane}}$ failures on average. The satellite failures in the operational plane follow a Poisson process. Thus, the times between consecutive orders from the operational plane follow an Erlang-$Q_{\text{plane}}$ distribution \cite{kim1991modeling}. The orders placed at all parking orbits combined are the superposition of the orders from all operational planes \cite{jakob2019optimal}. When $N_{\text{plane}}$ is sufficiently large (e.g., $\ge 20$), the superposition of those $N_{\text{plane}}$ Poisson processes is also considered a Poisson process with rate $N_{\text{plane}}\lambda_{\text{parking}}/Q_{\text{plane}}$ \cite{jakob2019optimal,ganeshan1999managing}. Therefore, as the parking orbits are equally distributed, it can be interpreted that each parking orbit is subject to a Poisson demand with rate $\lambda_{\text{parking}}$ due to its symmetricity. The demand rate of the parking orbit is derived below:
\begin{equation}\label{eqn: failure rate parking}
    \lambda_{\text{parking}} = N_{\text{plane}}   \frac{\lambda_{\text{plane}}}{Q_{\text{plane}}}   \frac{1}{N_{\text{parking}}}
\end{equation}
Again, in the dual-sourcing strategy, $Q_{\text{plane}}$ appears as an expected value, as derived in Eq. \eqref{eqn: expected reorder quantity plane}.

\subsubsection{Lead Time Distribution}
The depleted parking spares are directly supplied from the ground using a primary launch vehicle. Based on a similar approach as that used for in-plane spare replacement via direct injection, the lead time distribution for replenishing spares in the parking orbit ($l_3$) is modeled as a combination of an exponential distribution with a fixed term. This lead time accounts for the launch ordering process and the waiting period until the subsequent launch window. This study assumes that ground spares are always available, thus the availability factor is excluded from the lead time distribution modeling. The lead time for indirect injection ($w$) is represented as follows:
\begin{equation}\label{eqn: lead time parking}
    w \sim \exp\left(\mu_{\text{primary}}\right) + t_{\text{primary}}
\end{equation}
where $\mu_{\text{primary}}$ denotes the mean launch interval or the waiting time for the next launch, and $t_{\text{primary}}$ is the time required for processing the launch order.

\subsubsection{Expected Shortage}
The expected shortage of the parking spares can be derived using the basic technique introduced in Eq. \eqref{eqn: expected shortage expansion with Poisson}. Integrating the expected shortage over the lead time gives its expected value.
\begin{equation}\label{eqn: expected shortage parking}
    ES_{\text{parking}}=\int_{w\in D_{w}}ES\left(k_{R_1,3};\lambda_{\text{parking}},w\right)  l_3\left(w\right)dw
\end{equation}
Here, $w$ denotes the lead time modeled in Eq. \eqref{eqn: lead time parking}, $l_3$ and $D_w$ are the corresponding distribution and domain of $w$, respectively. Note that $ES_{\text{parking}}$ is derived in units of batches $Q_1$.

\subsubsection{Mean Stock Level}
The derivation of the mean stock level of the parking spare is similar to that of the expected shortage. Its relevant formulation is Eq. \eqref{eqn: mean stock level Poisson}, and integrating over lead time variable provides its expected value, as follows:
\begin{equation}\label{eqn: mean stock level parking}
    \overline{SL_{\text{parking}}}=\int_{w\in D_{l_3}} \overline{SL}\left(k_{R_1,3},k_{Q_1,3};\lambda_{\text{parking}}, w\right)  l_3\left(w\right)dw
\end{equation}
Similarly, $w$ represents the lead time modeled in Eq. \eqref{eqn: lead time parking}, with $l_3$ and $D_w$ being the associated distribution and domain of $w$, respectively. It should be noted that $\overline{SL_{\text{parking}}}$ is expressed in terms of units of batches $Q_1$ per time unit.

\subsection{Total Relevant Cost Model}
One of the primary objectives of the model developed in this subsection is to quantify the relevant costs associated with the spare strategy for maintaining the constellation. Following the techniques applied in the previous study \cite{jakob2019optimal} the total relevant cost is composed of four components: manufacturing ($c_{\text{manufacturing}}$), launch ($c_{\text{launch}}$), orbital maneuvering ($c_{\text{maneuvering}}$), and holding costs ($c_{\text{holding}}$). These components define the Total Expected Spare Strategy Annual Cost ($TESSAC$) \cite{jakob2019optimal}:
\begin{equation}\label{eqn: TESSAC}
TESSAC = c_{\text{manufacturing}} + c_{\text{launch}} + c_{\text{maneuvering}} + c_{\text{holding}}
\end{equation}

The remainder of this subsection models the cost component to construct the satellite mega-constellation optimal replenishment strategy.

\subsubsection{Satellite Manufacturing Cost}
The annual manufacturing cost is derived from the total number of plane failures over a year. As the Poisson process models satellite failure with a rate of $\lambda_{\text{plane}}$ for each orbital plane, the manufacturing cost $c_{\text{manufacturing}}$ is expressed as follows:
\begin{equation}\label{eqn: manufacturing cost}
c_{\text{manufacturing}} = c_{\text{sat}}  \lambda_{\text{plane}}  N_{\text{plane}}  N_t
\end{equation}
where $c_{\text{sat}}$ represents the manufacturing cost of a satellite, and $N_t$ is the number of time units in a year. The satellite failure rate in each orbital plane, $\lambda_{\text{plane}}$, is provided in Eq. \eqref{eqn: failure rate plane}.

\subsubsection{Launch Ordering Cost}
The annual launch cost is derived from the expected number of orders, specifically the orders from the parking orbits and orbital planes, through the indirect and direct channels, respectively. First, the annual number of launches of the primary vehicle dedicated for indirect injection to a parking orbit can be obtained by dividing the total annual demands by the reorder quantity in parking orbits:
\begin{equation}\label{eqn: annual order 3}
O_3 = \frac{\lambda_{\text{parking}} N_t}{k_{Q_1,3}}
\end{equation}

The annual number of launches of the auxiliary vehicle dedicated for direct injection to an orbital plane can be obtained similarly. Since auxiliary supply mode at in-plane spares triggers the direct supplement, the expected number of launches for the direct injection can be expressed as:

\begin{equation}\label{eqn: annual order 2}
O_2 = \frac{\lambda_{\text{plane}} N_t}{Q_{\text{plane}}} p_2 = \frac{\lambda_{\text{plane}} N_t}{Q_1 + Q_2 p_2} p_2
\end{equation}
Following these formulations, the total launch ordering cost is derived:
\begin{equation}
c_{\text{launch}} = c_{\text{auxiliary}} O_2 N_{\text{plane}} + c_{\text{primary}} O_3 N_{\text{parking}}
\end{equation}
Here, $c_{\text{primary}}$ and $c_{\text{auxiliary}}$ are the cost per launch of each channel. Since $O_3$ and $O_2$ are for a single parking orbit and a single orbital plane, respectively, $N_{\text{parking}}$ and $N_{\text{plane}}$ should be multiplied to obtain the total annual number of launches.

\subsubsection{In-plane Transfer Maneuvering Cost}
The annual maneuvering cost is relevant to the total fuel mass consumption to perform it for all orbital transfers required per year. To express the model, the annual number of orders from in-plane spares to the parking spares should be obtained. Similar to $O_3$ and $O_2$, it can be derived by dividing the expected annual demand at in-plane spare orbit by the expected reorder quantity:

\begin{equation}
    O_1 = \frac{\lambda_{\text{plane}}N_t}{Q_{\text{plane}}}
\end{equation}

The annual fuel requirement for orbital transfers is calculated by multiplying the number of orders ($O_1$) with the fuel consumption of a single maneuver ($m_{\text{fuel}}$). Using the cost conversion factor of the fuel mass ($\epsilon_{\text{fuel}}$) as proposed by \cite{jakob2019optimal}, the annual in-plane transfer maneuvering cost is derived as follows:
\begin{equation}
c_{\text{maneuvering}} = \epsilon_{\text{fuel}} O_1 N_{\text{plane}} Q_1 m_{\text{fuel}}
\end{equation}

\subsubsection{Spare Inventory Holding Cost}
The spare inventory holding cost is related to resources for holding the spare satellites in orbits, including orbital station-keeping maneuvers, monitoring the health of spare satellites, etc. Therefore, the annual holding cost for having spares in parking spare and in-plane spare is modeled by their mean stock levels with the parameter $h_s$, the annual holding cost per satellite:
\begin{equation}
    c_{\text{holding}} = h_s \left(\overline{SL_{\text{plane}}} N_{\text{plane}} + \overline{SL_{\text{parking}}} N_{\text{parking}} Q_1 \right)
\end{equation}
Since $\overline{SL_{\text{parking}}}$ is in units of batches $Q_1$, $Q_1$ multiplies the parking spare term.

The holding cost coefficient ($h_s$) relates to the actual cost of maintaining a spare satellite in orbit, but it can also represent a risk premium. Constellation operators may perceive keeping a spare satellite in orbit as risky due to environmental factors such as solar wind and orbital debris, which can potentially harm and cause malfunctions in their spares. Additionally, growing global concerns regarding space congestion can further burden constellation operators \cite{boley2021satellite}. These factors can be incorporated into the holding cost coefficient by adjusting its value to include both the actual holding cost and the risk component. 

\subsection{Model Validation}

\begin{table}[hbt!]
    \centering
    \caption{Common parameters of cases generated for validation study}
    \begin{tabular}{lccc}
        \hline\hline
        Parameter & Notation & Value & Unit\\\hline
        Fuel mass conversion coefficient & $\epsilon_{\text{fuel}}$ & 0.01 & M\$/kg\\
        Annual satellite holding cost & $h_{s}$ & 0.5 & M\$/satellite/year\\
        Satellite dry mass & $m_{\text{dry}}$ & 150 & kg\\
        Satellite production cost & $c_{\text{sat}}$ & 0.5 & M\$/satellite\\
        Exhaust velocity of satellite's propulsion system & $v_\text{ex}$ & 11.77 & km/s\\
        Mass flow rate of the satellite's propulsion system & $\dot{m}_{\text{low-thrust}}$ & 0.0013 & kg/s\\
        Primary vehicle launch cost & $c_{\text{normal}}$ & 67 & M\$\\
        Auxiliary vehicle launch cost & $c_{\text{auxiliary}}$ & 7.5 & M\$\\
        \hline
    \end{tabular}
    \label{tab: common parameters}
\end{table}

To validate the model, 25 cases were analyzed. The randomly sampled parameters for generating the test cases are summarized in Table \ref{tab: sampled trade space} with its ranges, and the fixed parameters that encompass the specifications of satellites with the cost-related parameters are outlined in Table \ref{tab: common parameters}. These parameters adhere to the model's construction assumptions, and the other detailed parameters can be found in Appendix D.

For the numerical experiments, the time window to initiate the auxiliary supply mode is parameterized as follows:
\begin{equation}\label{eqn: parameterized time window}
    t_w = t_{\text{trans}} + \alpha_w \frac{2\pi}{N_{\text{parking}}}\frac{1}{\left|\dot{\Omega}_{\text{relative}}\right|}
\end{equation}
where $t_{\text{trans}}$ is the in-plane transfer maneuver time, $N_{\text{parking}}$ is the number of parking orbits, and $\left|\dot{\Omega}_{\text{relative}}\right|$ is the difference in nodal precession between a parking orbit and an orbital plane of the constellation. With this parameterized time window, we can explain the concept of the time window intuitively. For instance, if $\alpha_w = j$ ($\le N_{\text{parking}}$), it indicates that the auxiliary supply mode can be activated until the spares from the $j^{\text{th}}$ closest parking orbit has arrived.

\begin{table}[hbt!]
    \centering
    \caption{Sampled trade space of varying parameters for validation study}
    \begin{tabular}{lccc}
        \hline\hline
        Parameter & Notation & Range & Unit\\\hline
        Launch order processing time of primary vehicle & $t_{\text{primary}}$ & $\left[4, 16\right]$ & weeks\\
        Mean waiting time to launch primary vehicle & $\mu_{\text{primary}}$ & $\left[4, 16\right]$ & weeks\\
        Launch order processing time of auxiliary vehicle & $t_{\text{auxiliary}}$ & $\left[1, 12\right]$ & weeks\\
        Mean waiting time to launch auxiliary vehicle & $\mu_{\text{auxiliary}}$ & $\left[1, 12\right]$ & weeks\\
        Orbital plane altitude & $h_{\text{plane}}$ & $\left[1000, 2000\right]$ & km\\
        Parking orbit altitude & $h_{\text{parking}}$ & $\left[400, 1000\right]$ & km\\
        Inclination & $i$ & $\left[40, 80\right]$ & deg\\
        Satellite failure rate & $\lambda_{\text{sat}}$ & $\left[0.05, 0.2\right]$ & failures/year\\
        Number of orbital planes & $N_{\text{plane}}$ & $\left[20, 40\right]$ & planes\\
        Number of parking orbits & $N_{\text{parking}}$ & $\left[1, 20\right]$ & planes\\
        Number of satellites per orbital plane & $N_{\text{sats}}$ & $\left[20, 60\right]$ & satellites/plane\\
        First reorder point at in-plane spare & $R_1$ & $\left[1, 20\right]$ & satellites\\
        Second reorder point at in-plane spare & $R_2$ & $\left[-2, 10\right]$ & satellites\\
        First order quantity at in-plane spare & $Q_1$ & $\left[1, 20\right]$ & satellites\\
        Second order quantity at in-plane spare & $Q_2$ & $\left[1, 10\right]$ & satellites\\
        Reorder point at parking spare & $k_{R_1,3}$ & $\left[1, 20\right]$ & $Q_{1}$\\
        Order quantity at parking spare & $k_{Q_1,3}$ & $\left[1, 20\right]$ & $Q_{1}$\\
        Time window parameter & $\alpha_{w}$ & $\left[0, 2\right]$ & -\\
        \hline\hline
    \end{tabular}
    \label{tab: sampled trade space}
\end{table}

To assess the model, a 30-year simulation was run on a weekly basis time step, and 100 Monte-Carlo simulations were conducted for each test case. The simulation outputs are compared to the results of the model, and the errors are calculated as following equations. The demand rate, mean stock level, and TESSAC are evaluated by the relative error through the following formula \cite{ganeshan1999managing}:

\begin{equation}\label{eqn:model value error}
\frac{\left|value_{\text{sim}} - value_{\text{model}}\right|}{value_{\text{sim}}} \times 100.
\end{equation}

The rational parameters, order fill rate and order cycle probability, are evaluated by the absolute error such as:
\begin{equation}\label{eqn:model rate error}
\left|value_{\text{sim}} - value_{\text{model}}\right|.
\end{equation}

\begin{table}[hbt!]
    \centering
    \caption{Averaged errors of result parameters}
    \begin{tabular}{lcc}
        \hline\hline
        Result Parameter & Notation & Mean Error \\\hline
        \textbf{Relative Errors (\%)} & & \\
        Demand rate at parking orbit & $\lambda_\text{sat}$ & 0.72\\
        Mean stock level at in-plane spare & $\overline{SL_\text{plane}}$ & 4.51\\
        Mean stock level at parking spare & $\overline{SL_\text{parking}}$ & 1.09\\
        Total relevant cost & $TESSAC$ & 1.45\\\addlinespace
        \textbf{Absolute Errors (\%p)} & & \\
        Order fill rate at in-plane spare & $\rho_\text{plane}$ & 0.04\\
        Order fill rate at parking spare & $\rho_\text{parking}$ & 0.10\\
        Two-order cycle probability at in-plane spare & $p_2$ &  0.49\\
        \hline\hline
    \end{tabular}
    \label{tab: averaged errors}
\end{table}

The results are summarized in Table \ref{tab: averaged errors}. These error percentages indicate that there are no significant differences between the analytic model and simulation values, as they are under 5\%. Therefore, the optimal replenishment spare policy can be established by employing the analytic expression according to the proposed framework, which approximates the expected behavior of the inventory model.

\section{Optimization Problem Formulation}
Two optimization problems are defined based on the developed inventory management model. In our context, only the satellite constellation replenishment strategies are concerned, excluding the constellation's design aspect. Thus, the problems are encoded given the specification of the constellation pattern and its satellites. This section concludes with the detailed context and formulation for each optimization problem.

\subsection{Problem 1: Optimization of Satellite Constellation Maintenance Strategy}
The satellite constellation operator is responsible for maintaining the service quality of their system. From this perspective, they can minimize the costs associated with maintaining the constellation system's performance criteria. The service provider can adjust the configuration of the parking orbits for spare satellites and the replenishment parameters for each spare in the preliminary design phase. By regulating these factors, the optimal balance among various cost components can be identified to minimize the total relevant cost while preserving service quality. The cost components for the satellite constellation replenishment problem are considered as; launch ordering costs, orbit transfer maneuver costs, spare holding costs, and the satellite manufacturing cost. The problem is abstracted as the following optimization model ($\mathbf{P}_{\text{OR}}$).
\newline\newline\noindent
($\mathbf{P}_{\text{OR}}$) Optimal Replenishment Strategy with Dual Supply Modes

\begin{equation}\label{eqn: P_OR obj}
    \underset{\mathbf{R},\mathbf{Q},\mathbf{x}_{\text{parking}},\alpha_w}{\min}\ TESSAC\left(\mathbf{R},\mathbf{Q}, \mathbf{x}_{\text{parking}};\mathbf{q}\right)
\end{equation}
\newline\noindent subject to
\begin{equation}\label{eqn: P_OR const 1}
    \rho_{\text{parking}}\left(\mathbf{R},\mathbf{Q},\mathbf{x}_{\text{parking}};\mathbf{q}\right)\ge \rho_{\text{parking}}^{\text{req}}
\end{equation}
\begin{equation}\label{eqn: P_OR const 2}
    \rho_{\text{plane}}\left(\mathbf{R},\mathbf{Q},\mathbf{x}_{\text{parking}};\mathbf{q}\right)\ge \rho_{\text{plane}}^{\text{req}}
\end{equation}
\begin{equation}\label{eqn: P_OR const 3}
    R_1-R_2 \ge 1
\end{equation}
\begin{equation}\label{eqn: P_OR const 4}
    Q_1  k_{Q_1,3} \le Q_{3,\text{max}}
\end{equation}
\begin{equation}\label{eqn: P_OR const 5}
    Q_2 \le Q_{2,\text{max}}
\end{equation}
\begin{equation}\label{eqn: P_OR const 6}
    R_1 \le Q_{\text{plane}}
\end{equation}
\begin{equation}\label{eqn: P_OR const 7}
    k_{R_1,3} \le k_{Q_1,3}
\end{equation}
\begin{equation}\label{eqn: P_OR const 8}
    \mathbf{R} =\left[R_1,R_2,k_{R_1,3}\right]\in \mathbb{Z}^3
\end{equation}
\begin{equation}\label{eqn: P_OR const 9}
    \mathbf{Q} =\left[Q_1,Q_2,k_{Q_1,3}\right]\in \mathbb{Z}_{++}^3
\end{equation}
\begin{equation}\label{eqn: P_OR const 10}
    \mathbf{x}_{\text{parking}}=\left[N_{\text{parking}},h_{\text{parking}}\right]\in \mathcal{N}_{\text{parking}}\times \mathcal{H}_{\text{parking}}
\end{equation}
\begin{equation}\label{eqn: P_OR const 11}
    \alpha_w \in \mathbb{R}_{+}
\end{equation}
Here, $\mathbf{R}$, $\mathbf{Q}$, and $\mathbf{x}_{\text{parking}}$ are decision variables related to the reorder point, the order quantity, parking orbits configuration, and $\mathbf{q}$ is the vector of the parameters related to the corresponding function, respectively. As explained above, the constellation pattern, specification of the satellite, and its relevant cost coefficients are parameters. The objective function of the optimization problem is to minimize $TESSAC$.

In addition, constraints for maintaining the service levels of parking and in-plane spares are given in Eqs. \eqref{eqn: P_OR const 1} and \eqref{eqn: P_OR const 2}, respectively. Eq. \eqref{eqn: P_OR const 3} ensures the constraint on reorder points of normal and auxiliary supply modes, while Eqs. \eqref{eqn: P_OR const 4} and \eqref{eqn: P_OR const 5} enforce the total launch quantity within the launch capacity limit of the primary vehicle ($Q_{3,\text{max}}$) and auxiliary vehicle ($Q_{2,\text{max}}$). Eqs. \eqref{eqn: P_OR const 6} and \eqref{eqn: P_OR const 7} are introduced as supplementary constraints for maintaining the cyclic behavior of stocks. Finally, Eqs. \eqref{eqn: P_OR const 8} to \eqref{eqn: P_OR const 11} define the decision variables, with $\mathcal{N}_{\text{parking}}$ and $\mathcal{H}_{\text{parking}}$ representing the feasible sets for parking configuration variables.

\subsection{Problem 2: Valuation of the Direct Injection - Auxiliary Launch Option}
In the proposed injection strategy, the responsive launch option (direct injection via auxiliary launcher) can provide lower total maintenance cost and flexibility in the satellite constellation replenishment design. To evaluate the option, we assume that direct injection is only served by a small responsive auxiliary launcher and formulated the problem from the customer's perspective.

The value of an auxiliary launch option can be quantified by the \textit{willingness to pay} of the constellation operator. This represents the highest amount a constellation operator is willing to pay for the direct launch option when the auxiliary launcher can be utilized for the constellation maintenance. The optimization problem $\mathbf{P}_{\text{VA}}$ is introduced to encode this context.
\newline\newline\noindent
($\mathbf{P}_{\text{VA}}$) Value of auxiliary launch option
\begin{equation}\label{eqn: P_VA obj}
    \underset{\mathbf{R},\mathbf{Q},\mathbf{x}_{\text{parking}},\alpha_w,c_{\text{auxiliary}}}{\max}\ c_{\text{auxiliary}}
\end{equation}

\noindent subject to
\begin{equation}\label{eqn: P_VA const 1}
    TESSAC\left(\mathbf{R},\mathbf{Q},\mathbf{x}_{\text{parking}},c_{\text{auxiliary}};\mathbf{q}\right)\le TESSAC_{\text{ref}}
\end{equation}
\begin{equation}\label{eqn: P_VA const 2}
    k_{Q_1,3}p_2\left(\mathbf{R},\mathbf{Q},\mathbf{x}_{\text{parking}};\mathbf{q}\right) \ge \eta
\end{equation}
\begin{equation}\label{eqn: P_VA const 3}
    c_{\text{auxiliary}}\in \mathbb{R}_{++}
\end{equation}
\begin{center}
    Eqs. \eqref{eqn: P_OR const 1}-\eqref{eqn: P_OR const 11}
\end{center}
In this formulation, $c_{\text{auxiliary}}$ is an additional decision variable to be determined. The objective function is the maximization of the auxiliary vehicle launch cost, which is the willingness to pay for the direct injection. The constraints in $\mathbf{P}_{\text{SO}}$ are also considered. In addition to these constraints, Eq. \eqref{eqn: P_VA const 1} restricts the budget for maintenance relative to the reference annual maintenance cost ($TESSAC_{\text{ref}}$). $TESSAC_{\text{ref}}$ can be set as the constellation operator's budget or from other maintenance strategies. Eq. \eqref{eqn: P_VA const 2} establishes a constraint that ensures the relative usage of auxiliary launcher usage to primary launcher usage, which is represented by the multiplication of the number of in-plane stock cycles enabled by a primary launcher and the probability of triggering a two-order cycle ($k_{Q_1,3}p_2$), meets or exceeds a specified minimum threshold ($\eta$).

To summarize, by solving this problem, the value of the direct injection can be determined through the constellation operator's willingness to pay, under the reference budget and the targeted sales share of the launch service market.

\section{Case Study}
This section summarizes case studies for the proposed optimization problems defined in the previous section. The case studies include the worked examples of the proposed problems with the real-world sized constellation and its satellite specifications. The results can offer the baseline for designing and analyzing the satellite constellation maintenance strategy to the relevant stakeholders.

\subsection{Case Study 1: \texorpdfstring{$\left[\mathbf{P}_{\text{OR}}\right]$}{} Optimization of Satellite Constellation Maintenance Strategy}
Case Study 1 is the worked example for $\mathbf{P}_{\text{OR}}$. This conducts two launch injection options: an indirect injection via a primary vehicle and a direct injection through an auxiliary vehicle. The primary vehicle is characterized by superior payload capacity with its lower specific cost (ratio of the launch cost to its maximum payload capacity) but an extended lead time, and an auxiliary vehicle is characterized by exhibiting inverse characteristics. Table \ref{tab: fixed case study parameters 1} summarizes the parameters for Case Study 1. Note that the $\alpha_w$ and $Q_2$ are set as fixed parameters, which are decision variables in the general formulation $\mathbf{P}_{\text{OR}}$.

\begin{table}[hbt!]
    \centering
    \caption{Fixed case study parameters of Case Study 1}
    \begin{tabular}{lccc}
        \hline\hline
        Parameter & Notation & Value & Unit\\\hline
        Fuel mass conversion coefficient & $\epsilon_{\text{fuel}}$ & 0.01 & M\$/kg\\
        Satellite dry mass & $m_{\text{dry}}$ & 150 & kg\\
        Satellite production cost & $c_{\text{sat}}$ & 0.5 & M\$/satellite\\
        Exhaust velocity of satellite's propulsion system & $v_{ex}$ & 11.77 & km/s\\
        Mass flow rate of the satellite's propulsion system & $\dot{m}_{\text{low-thrust}}$ & 0.0013 & kg/s\\
        Annual satellite holding cost & $h_s$ & 0.5 & M\$/satellite/year\\
        Primary vehicle launch cost & $c_{\text{primary}}$ & 67 & M\$\\
        Launch order processing time of primary vehicle & $t_{\text{primary}}$ & 12 & weeks\\
        Mean waiting time to launch primary vehicle & $\mu_{\text{primary}}$ &  8  & weeks\\
        Primary vehicle launch capacity & $Q_{3,\text{max}}$ & 40 & satellites\\
        Auxiliary vehicle launch cost & $c_{\text{primary}}$ & 7.5 & M\$\\
        Launch order processing time of auxiliary vehicle & $t_{\text{auxiliary}}$ & 2 & weeks\\
        Mean waiting time to launch auxiliary vehicle & $\mu_{\text{auxiliary}}$ & 2 & weeks\\
        Orbital plane altitude & $h_{\text{plane}}$ & 1200 & km\\
        Inclination & $i$ & 60 & deg\\
        Satellite failure rate & $\lambda_{\text{sat}}$ & 0.2 & failures/year\\
        Number of orbital planes & $N_{\text{plane}}$ & 40 & planes\\
        Number of satellites per orbital plane & $N_{\text{sats}}$ & 40 & satellites/plane\\
        Required order fill rate at in-plane spare & $\rho_{\text{plane}}^{\text{req}}$ & 0.98 & -\\
        Required order fill rate at parking spare & $\rho_{\text{parking}}^{\text{req}}$ & 0.98 & -\\
        Time window parameter & $\alpha_w$ & 1 & - \\
        Second order quantity at in-plane spare & $Q_2\left(=Q_{2,\text{max}}\right)$ & 2 & satellites\\
        \hline\hline
    \end{tabular}
    \label{tab: fixed case study parameters 1}
\end{table}
The decision variables comprise the problem $\mathbf{P}_{\text{OR}}$ is illustrated in Table \ref{tab: decision variables 1} with its ranges. The altitude of the parking orbit $\mathcal{H}_{\text{parking}}$ is quantized by 50 km intervals. So the possible values of the parking orbit's altitude are converted into the integers within the given ranges:
\begin{equation}\label{eqn: possible values of parking orbit altitude}
\mathcal{H}_{\text{parking}} = \left\{400, 450, \dots, 1100\right\}
\end{equation}

\begin{table}[hbt!]
    \centering
    \caption{Decision variables for optimization of Case Study 1}
    \begin{tabular}{lccc}
        \hline\hline
        Decision Variable & Notation & Range & Unit\\\hline
        Number of parking orbits & $N_{\text{parking}}$ & $\left[1, 20\right]$ & planes\\
        Parking orbit altitude & $h_{\text{parking}}$ & $\left[400, 1100\right]$ & km\\
        First reorder point at in-plane spare & $R_1$ & $\left[1, 20\right]$ & satellites\\
        Second reorder point at in-plane spare & $R_2$ & $\left[-2, 10\right]$ & satellites\\
        First order quantity at in-plane spare & $Q_1$ & $\left[1, 40\right]$ & satellites\\
        Reorder point at parking spare & $k_{R_1,3}$ & $\left[1, 20\right]$ & $Q_{1}$\\
        Order quantity at parking spare & $k_{Q_1,3}$ & $\left[1, 40\right]$ & $Q_{1}$\\
        \hline\hline
    \end{tabular}
    \label{tab: decision variables 1}
\end{table}

\begin{table}[hbt!]
    \centering
    \caption{Optimal solution and selected result parameters of Case Study 1}
    \begin{tabular}{lc}
        \hline\hline
         & Value \\\hline
        \textbf{Optimal Solution} & \\
        $N_{\text{parking}}$ (planes) & 9 \\
        $h_{\text{parking}}$ (km) & 650 \\
        $R_1$ (satellites) & 3 \\
        $R_2$ (satellites) & -2 \\
        $Q_1$ (satellites) & 5 \\
        $k_{R_1,3}$ ($Q_1$) & 5 \\
        $k_{Q_1,3}$ ($Q_1$) & 8 \\\addlinespace
        \textbf{Result Parameter} & \\
        $\overline{SL_{\text{plane}}}$ (satellites) & 4.9 \\
        $\overline{SL_{\text{parking}}}$ ($Q_1$) & 6.8 \\
        $TESSAC$ (M\$/year) & 962.1 \\
        $\rho_{\text{plane}}$ (\%) & 98.2 \\
        $\rho_{\text{parking}}$ (\%) & 98.1 \\
        $p_2$ (\%) & 1.2 \\
        \hline\hline
    \end{tabular}
    \label{tab: optimal solution 1}
\end{table}

The Genetic Algorithm (GA) is employed to solve the generated instances, and the results are presented in Table \ref{tab: optimal solution 1}. The selection of $R_2 = -2$, which is the lower bound, and its difference of 5 from $R_1$, results in low $p_2$ and $k_{Q_1,3}p_2$ values. Specifically, the value of $k_{Q_1,3}p_2$ about 0.1 implies that the auxiliary launcher is used only about once for every 10 uses of the primary launcher. Consequently, in this scenario, the dominant strategy becomes utilization of the primary launcher and indirectly supplying spares by transferring them from parking orbits to in-plane positions. In other words, the direct channel provided by the auxiliary launcher becomes a significantly less favored option.

\subsection{Case Study 2: \texorpdfstring{$\left[\mathbf{P}_{\text{OR}}\right]$}{} Parametric Study with Varying Holding Cost}
The optimal replenishment strategy depends mainly on the cost structure. The indirect channel-dominant strategy relies more heavily on the primary launcher, which is the solution of the previous case study. It is characterized by lower specific launch costs but higher expenses for holding spares in parking orbits and transferring them to in-plane positions. This strategy becomes preferable when the cost savings from using the more efficient primary launcher outweigh the additional expenses of storage and transfer. Conversely, the direct channel utilizes the auxiliary launcher more frequently. It involves higher specific launch costs but minimizes the expenses associated with holding stocks in parking orbits and transferring them. This approach becomes dominant when the higher launch costs are offset by savings in storage and transfer expenses.

Case Study 2 explores the trade-off between indirect and direct replenishment strategies by examining five instances with varying holding cost coefficients. This approach allows for the analysis of optimal solutions across diverse contextual situations and perceived risks associated with holding spare satellites in orbit. The study incrementally increases the holding cost from the reference value established in Case Study 1 (instance \#0), raising it by 0.1M\$/satellite/year, up to a maximum increase of 1.0M\$/satellite/year. The instances are distinguished from \#1 to \#5. To streamline the analysis and clearly observe trends, the study focuses solely on variations in the replenishment policy for spares by adopting the parking orbit configuration from the optimal solution from Case Study 1. The decision variables for this case study and the resulting optimal solutions are presented in Table \ref{tab: decision variables 2} and Table \ref{tab: optimal solution 2}, respectively.

\begin{table}[hbt!]
    \centering
    \caption{Decision variables for optimization of Case Study 2}
    \begin{tabular}{lccc}
        \hline\hline
        Decision Variable & Notation & Range & Unit\\\hline
        First reorder point at in-plane spare & $R_1$ & $\left[1, 20\right]$ & satellites\\
        Second reorder point at in-plane spare & $R_2$ & $\left[-2, 10\right]$ & satellites\\
        First order quantity at in-plane spare & $Q_1$ & $\left[1, 40\right]$ & satellites\\
        Reorder point at parking spare & $k_{R_1,3}$ & $\left[1, 20\right]$ & $Q_{1}$\\
        Order quantity at parking spare & $k_{Q_1,3}$ & $\left[1, 40\right]$ & $Q_{1}$\\
        \hline\hline
    \end{tabular}
    \label{tab: decision variables 2}
\end{table}

\begin{table}[hbt!]
    \centering
    \caption{Optimal solution and selected result parameters of Case Study 2}
    \begin{tabular}{lcc}
        \hline\hline
         & \multicolumn{2}{c}{Value for Instance}\\
         & \#0, \#1, \#2, \#3 & \#4, \#5\\\hline
        \textbf{Optimal Solution} & &\\
        $R_1$ (satellites) & 3 & 3 \\
        $R_2$ (satellites) & -2 & 0 \\
        $Q_1$ (satellites) & 5 & 3 \\
        $k_{R_1,3}$ ($Q_1$) & 5 & 7 \\
        $k_{Q_1,3}$ ($Q_1$) & 8 & 13 \\\addlinespace
        \textbf{Result Parameter} & & \\
        $\overline{SL_{\text{plane}}}$ (satellites) & 4.9 & 4.1\\
        $\overline{SL_{\text{parking}}}$ ($Q_1$) & 6.8 & 9.9\\
        $\rho_{\text{plane}}$ (\%) & 98.2 & 98.2\\
        $\rho_{\text{parking}}$ (\%) & 98.1 & 98.3\\
        $p_2$ (\%) & 1.2 & 13.6\\
        \hline\hline
    \end{tabular}
    \label{tab: optimal solution 2}
\end{table}
The output shows a significant shift in the solution between holding costs of 0.8 and 0.9. Instances \#0 through \#3 maintain the solution derived in Case Study 1, while instances \#4 and \#5 exhibit a new but identical solution. This abrupt transition is attributed to the discrete nature of the optimization problem under consideration. Specifically, even though the reorder point ($R_1$), remains the same for all instances, the reorder point of the auxiliary supply mode ($R_2$) increases from -2 to 0, and the order quantity of the normal supply mode ($Q_1$) decreases from 5 to 3, implying that the policy changes to one that relies more on the auxiliary supply mode.

This replenishment policy transition derives significant changes in the resultant inventory parameters. The total mean stock level---computed as the sum of the total mean stock level in in-plane spare ($\overline{SL_{\text{plane}}}N_{\text{plane}}$) and that in parking spare ($\overline{SL_{\text{parking}}}N_{\text{parking}}Q_1$)---decreases from 503.0 to 431.9, about 14.1\%, and the relative usage of auxiliary launchers per primary launch ($k_{Q_1,3}$$p_2$) increases from 0.1 to 1.8, about 18 times.

\begin{figure}[hbt!]
    \centering
    \includegraphics[width=0.8\linewidth, page=10]{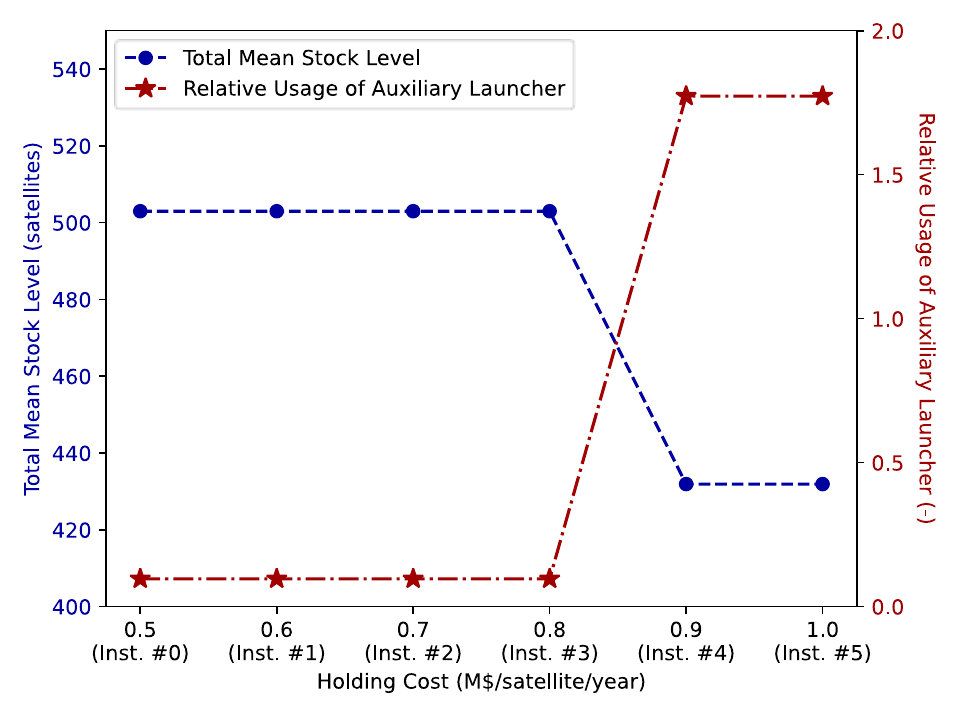}
    \caption{Total mean stock level and relative usage of the auxiliary launcher of the optimal solutions in Case Study 2}
    \label{fig: case study 2 result}
\end{figure}

\subsection{Case Study 3: \texorpdfstring{$\left[\mathbf{P}_{\text{VA}}\right]$}{} Valuation of the Direct Injection - Auxiliary Launch Option}

The third worked example focuses on evaluating the direct channel from the perspective of the auxiliary launch service. The scenarios follow the parameters in the previous case study with instance \#0 and \#5. Baseline strategies with a single-channel supply chain are established first, to assess the value of their service. The optimal replenishment strategy with the single-channel model is determined and presented in Table \ref{tab: optimal solution single}. This optimization problem is formulated using $\mathbf{P}_{\text{OR}}$ without considering the auxiliary supply mode. The single-channel inventory control model used in this analysis is adopted from the previous research \cite{jakob2019optimal}. The obtained optimal solution ($TESSAC$ and parking orbit configuration - $N_{\text{parking}}$, $h_{\text{parking}}$) are used as the $TESSAC{\text{ref}}$ and the parking orbit configuration of $\mathbf{P}_{\text{VA}}$, respectively.

\begin{table}[hbt!]
    \centering
    \caption{Optimal solution and selected result parameters of single-channel strategy for Case Study 3}
    \begin{tabular}{lcc}
        \hline\hline
         & \multicolumn{2}{c}{Value for Instance}\\
         & \#0 & \#5 \\\hline
        \textbf{Optimal Solution} & & \\
        $N_{\text{parking}}$ & 8 & 11 \\
        $h_{\text{parking}}$ (km) & 650 & 550 \\
        $R_1$ (satellites) & 3 & 3 \\
        $Q_1$ (satellites) & 4 & 4 \\
        $k_{R_1,3}$ ($Q_1$) & 7 & 5 \\
        $k_{Q_1,3}$ ($Q_1$) & 10 & 10 \\\addlinespace
        \textbf{Result Parameter} & & \\
        $\rho_{\text{plane}}$ (\%) & 98.1 & 98.6 \\
        $\rho_{\text{parking}}$ (\%) & 98.4 & 98.3\\
        $TESSAC$ (M\$/year) & 931.0 & 1221.7 \\
        \hline\hline
    \end{tabular}
    \label{tab: optimal solution single}
\end{table}

In addition to the reference budget and the parking orbit configuration, the required market share is set as $\eta=2$. This value implies that at least two auxiliary launch services are expected per primary launch service. With the parameters and the decision variables as defined in Table \ref{tab: decision variables 3}, the optimal solution is obtained and presented in Table \ref{tab: optimal solution 3}.

\begin{table}[hbt!]
    \centering
    \caption{Decision variables for optimization of Case Study 3}
    \begin{tabular}{lccc}
        \hline\hline
        Variables & Notation & Range & Unit\\\hline
        First reorder point at in-plane spare & $R_1$ & $\left[1, 20\right]$ & satellites\\
        Second reorder point at in-plane spare & $R_2$ & $\left[-2, 20\right]$ & satellites\\
        First order quantity at in-plane spare & $Q_1$ & $\left[1, 40\right]$ & satellites\\
        Second order quantity at in-plane spare & $Q_2$ & $\left[1, 2\right]$ & satellites\\
        Reorder point at parking spare & $k_{R_1,3}$ & $\left[1, 20\right]$ & $Q_{1}$\\
        Order quantity at parking spare & $k_{Q_1,3}$ & $\left[1, 40\right]$ & $Q_{1}$\\
        Auxiliary vehicle launch cost & $c_{\text{auxiliary}}$ & $\left[0, 5\right]$ & M\$\\
        \hline\hline
    \end{tabular}
    \label{tab: decision variables 3}
\end{table}

The result reveals that the constellation operator's willingness to pay for an auxiliary launch service as a direct supply channel ($c_{\text{auxiliary}}^*$) are 3.7M\$ and 4.3M\$, which are 49\% and 57\% of its initial price 7.5M\$, respectively. This indicates that the price of the auxiliary launch service should be reduced to achieve the desired relative usage of their launch service. However, it is noteworthy that $k_{Q_1,3}p_2$ are over 9 in both instances, implying overachievements of the target value of 2, which can be attributed to the discrete nature of the problem.

In another aspect, when examining the specific launch costs, the primary launch vehicle demonstrates a cost of 1.675 M\$/satellite (=$c_{\text{primary}}/Q_{3,\text{max}}$), whereas the auxiliary launch vehicle exhibits costs of 1.85M\$/satellite and 2.15 M\$/satellite (=$c_{\text{auxiliary}}^*/Q_{2,\text{max}}$). This comparison reveals that the auxiliary launcher's specific cost is slightly higher than that of the primary launcher, and the gap increases as the holding cost increases.

\begin{table}[hbt!]
    \centering
    \caption{Optimal solutions and selected inventory parameters of Case Study 3}
    \begin{tabular}{lcc}
        \hline\hline
         & \multicolumn{2}{c}{Value for Instance}\\
         & \#0 & \#5 \\\hline
        \textbf{Optimal Solution} & & \\
        $R_1$ (satellites) & 3 & 2\\
        $R_2$ (satellites) & 2 & 1\\
        $Q_1$ (satellites) & 2 & 3\\
        $k_{R_1,3}$ ($Q_1$) & 7 & 4\\
        $k_{Q_1,3}$ ($Q_1$) & 20 & 13\\
        $c_{\text{auxiliary}}$ (M\$) & 3.7 & 4.3 \\\addlinespace
        \textbf{Result Parameter} & & \\
        $\overline{SL_{\text{plane}}}$ (satellites) & 4.9 & 4.4\\
        $\overline{SL_{\text{parking}}}$ ($Q_1$) & 13.1 & 8.5\\
        $TESSAC$ (M\$/year) & 929.0 & 1218.4 \\
        $\rho_{\text{plane}}$ (\%) & 99.5 & 98.3\\
        $\rho_{\text{parking}}$ (\%) & 98.5 &  98.1\\
        $p_2$ (\%) & 71.4 & 70.3\\
        \hline\hline
    \end{tabular}
    \label{tab: optimal solution 3}
\end{table}

\section{Conclusion}
This study introduces an inventory management model for spare satellite strategy in mega-constellations, incorporating dual supply modes: normal and auxiliary. The framework combines an indirect channel for normal supply and a direct channel for auxiliary supply for optimizing the spare satellites replenishment policy. It provides analytical models for various inventory-related parameters, including mean stock level, expected shortage, and order cycle probabilities, followed by a cost model. To define the inventory management model, parametric replenishment policies, $(s,Q)$ and $(R_1,R_2,Q_1,Q_2)$ with a time window, are adopted and modified for this specific context of the problem.

Based on the proposed model, two optimization problems are proposed and formulated from different perspectives: the constellation operator and the launch service provider. Case studies with detailed contexts follow the formulation, and the worked examples are conducted in each study. The results align with qualitative intuitions; for example, increased holding costs push the constellation operator to rely more on a direct injection strategy than maintaining spare satellites in orbit via indirect injection. Furthermore, the results provide specific solution values that could offer a baseline in the decision-making processes of stakeholders.

Incorporating environmental factors into the model and optimization framework presents a promising avenue for future research. As space becomes increasingly active, it is becoming more congested and potentially hostile to orbital assets. Future studies could focus on integrating environmental considerations in several ways: incorporating environmental factors into cost elements, selecting parking orbits by accounting for environmental factors, or proposing new supply chain models that utilize fewer orbital resources. These approaches would address the growing concerns of space sustainability while enhancing the model's relevance in the evolving context of space operations.

\section*{Appendix}
To simplify the mathematical expression, $f\left(x;\tau\right)$ substitutes $f_{\lambda_{\text{plane}}}\left(x;\tau\right)$ by omitting the subscript $\lambda_{\text{plane}}$ from Appendix A to C.

\subsection{Expected Shortage}
\begin{align}
    &ES^1\left(y;R_1,R_2,Q_1,Q_2,t_w\right)\nonumber\\
    &=\begin{cases}\int_{R_1^+}^{\max\left(R_1-R_2,R_1^+\right)}\left(x-R_1\right)  f\left(x;y\right)dx, & y \le t_w\\
    \int_{0}^{R_1-R_2}\int_{\left[R_1-x\right]^+}^{\infty}\left(x'-R_1+x\right)  f\left(x';y-t_w\right)  f\left(x;t_w\right) dx'dx, & y> t_w
    \end{cases}
\end{align}
\begin{align}
    &ES^{1,2}\left(y,t,z;R_1,R_2,Q_1,Q_2,t_w\right)\nonumber\\
    &= \int_{R_2^+}^{\infty}\left(x-R_2\right)  f\left(x;y-t\right) dx\nonumber\\
    &\ \ \ +\int_{\left[R_2+Q_1\right]^+}^{\infty}\int_{0}^{\infty} x'  f\left(x';z+t-y\right)  f\left(x;y-t\right) dx'dx\nonumber\\
    &\ \ \ +\int_{0}^{\left[R_2+Q_1\right]^+}\int_{\left[R_2+Q_1-x\right]^+}^{\infty}\left(x'-R_2-Q_1+x\right)  f\left(x';z+t-y\right)  f\left(x;y-t\right) dx'dx
\end{align}
\begin{align}
    &ES^{2,1}\left(y,t,z;R_1,R_2,Q_1,Q_2,t_w\right)\nonumber\\
    &= \int_{R_2^+}^{\infty}\left(x-R_2\right)  f\left(x;z\right) dx\nonumber\\
    &\ \ \ +\int_{\left[R_2+Q_1\right]^+}^{\infty}\int_{0}^{\infty} x'  f\left(x';y-z-t\right)  f\left(x;z\right) dx'dx\nonumber\\
    &\ \ \ +\int_{0}^{\left[R_2+Q_2\right]^+}\int_{\left[R_2+Q_2-x\right]^+}^{\infty}\left(x'-R_2-Q_2+x\right)  f\left(x';y-z-t\right)  f\left(x;z\right) dx'dx
\end{align}

\subsection{Average Cycle Stock}
\begin{align}
    &CS^1\left(y;R_1,R_2,Q_1,Q_2,t_w\right)\nonumber\\
    &= \begin{cases}\int_{0}^{R_1-R_2}y  \left(R_1-\frac{x}{2}\right)  f\left(x;y\right)dx + \int_{0}^{R_1-R_2}\frac{Q_1-x}{\lambda_{\text{plane}}}  \left(R_1+\frac{Q_1-x}{2}\right)  f\left(x;y\right)dx, & y\le t_w\\
    \int_{0}^{R_1-R_2}t_w  \left(R_1-\frac{x}{2}\right)  f\left(x;t_w\right)dx\\
    \ \ \ \ \ +\int_{0}^{R_1-R_2}\int_{0}^{\infty}\left(y-t_w\right)  \left(R_1-x-\frac{x'}{2}\right)  f\left(x';y-t_w\right)dx'dx\\
    \ \ \ \ \ +\int_{0}^{R_1-R_2}\int_{0}^{\infty}\frac{Q_1-x-x'}{\lambda_{\text{plane}}}  \left(R_1+\frac{Q_1-x-x'}{2}\right)  f\left(x';y-t_w\right)dx'dx, & y > t_w
    \end{cases}
\end{align}
\begin{align}
    &CS^{1,2}\left(y,t,z;R_1,R_2,Q_1,Q_2,t_w\right)\nonumber\\
    &= t  \frac{R_1+R_2}{2}+\nonumber\\
    &\ \ \ \int_{0}^{\infty}\left(y-t\right)  \left(R_2-\frac{x}{2}\right)  f\left(x;y-t\right)dx\nonumber\\
    &\ \ \ +\int_{0}^{\infty}\int_{0}^{\infty}\left(t+z-y\right)  \left(R_2+Q_1-x-\frac{x'}{2}\right)  f\left(x';t+z-y\right)  f\left(x;y-t\right)dx'dx\nonumber\\
    &\ \ \ +\int_{0}^{\infty}\int_{0}^{\infty}\frac{R_2+Q_1+Q_2-x-x'-R_1}{\lambda_{\text{plane}}} \frac{R_2+Q_1+Q_2-x-x'+R_1}{2}  f\left(x';t+z-y\right)  f\left(x;y-t\right)dx'dx\label{eqn: two-order average cycle stock 1}
\end{align}
\begin{align}
    &CS^{2,1}\left(y,t,z;R_1,R_2,Q_1,Q_2,t_w\right)\nonumber\\
    &=t  \frac{R_1+R_2}{2}+\nonumber\\
    &\ \ \ \int_{0}^{\infty}z  \left(R_2-\frac{x}{2}\right)  f\left(x;z\right)dx\nonumber\\
    &\ \ \ +\int_{0}^{\infty}\int_{0}^{\infty}\left(y-t-z\right)  \left(R_2+Q_2-x-\frac{x'}{2}\right)  f\left(x';y-t-z\right)  f\left(x;z\right)dx'dx\nonumber\\
    &\ \ \ +\int_{0}^{\infty}\int_{0}^{\infty}\frac{R_2+Q_1+Q_2-x-x'-R_1}{\lambda_{\text{plane}}}  \frac{R_2+Q_1+Q_2-x-x'-R_1}{2}  f\left(x';y-t-z\right)  f\left(x;z\right) dx'dx\label{eqn: two-order average cycle stock 2}
\end{align}

\subsection{Order Cycle Probability}
\begin{equation}\label{eqn: one-order cycle probability}
    p_1 = \int_{0}^{\infty}\int_{0}^{R_1-R_2}f\left(x;\min\left(y,t_w\right)\right)  l_1\left(y\right)dxdy
\end{equation}

\begin{equation}\label{eqn: two-order cycle probability}
    p_2 = \int_{0}^{\infty}\int_{R_1-R_2}^{\infty}f\left(x;\min\left(y,t_w\right)\right)  l_1\left(y\right)dxdy
\end{equation}

\subsection{Parameters for validation study}
\begin{table}[hbt!]
    \centering
    \caption{Parameters for validation study -- 1}
    \begin{tabular}{lrrrrrrrr}
        \hline\hline
        Inst. \# & $t_{\text{primary}}$ & $\mu_{\text{primary}}$ & $t_{\text{auxiliary}}$ & $\mu_{\text{auxiliary}}$ & $h_{\text{plane}}$ & $h_{\text{parking}}$ & $i$ & $\lambda_{\text{sat}}$\\\hline
        1   & 12    & 6     & 8     & 9     & 1600  & 650   & 80    & 0.08\\
        2   & 8     & 14    & 4     & 4     & 1000  & 600   & 56    & 0.06\\
        3   & 4     & 12    & 8     & 11    & 1450  & 750   & 70    & 0.05\\
        4   & 16    & 8     & 12    & 3     & 1200  & 550   & 74    & 0.11\\
        5   & 7     & 16    & 6     & 7     & 1600  & 1000  & 40    & 0.15\\
        6   & 14    & 10    & 3     & 9     & 1100  & 450   & 64    & 0.2\\
        7   & 4     & 15    & 2     & 12    & 1650  & 600   & 48    & 0.13\\
        8   & 6     & 12    & 10    & 8     & 1450  & 850   & 48    & 0.19\\
        9   & 14    & 7     & 1     & 10    & 1250  & 400   & 36    & 0.18\\
        10  & 11    & 10    & 11    & 4     & 1300  & 600   & 78    & 0.07\\
        11  & 5     & 11    & 6     & 6     & 1050  & 950   & 42    & 0.09\\
        12  & 10    & 13    & 2     & 5     & 2000  & 700   & 68    & 0.17\\
        13  & 7     & 5     & 5     & 11    & 1350  & 900   & 54    & 0.18\\
        14  & 15    & 4     & 11    & 7     & 1950  & 1000  & 60    & 0.17\\
        15  & 10    & 6     & 7     & 12    & 1150  & 700   & 62    & 0.06\\
        16  & 8     & 9     & 4     & 6     & 1850  & 650   & 46    & 0.14\\
        17  & 9     & 14    & 9     & 1     & 1400  & 850   & 72    & 0.14\\
        18  & 15    & 4     & 12    & 5     & 1300  & 900   & 36    & 0.08\\
        19  & 13    & 5     & 1     & 1     & 1250  & 500   & 34    & 0.05\\
        20  & 9     & 11    & 9     & 10    & 1900  & 950   & 58    & 0.1\\
        21  & 5     & 9     & 1     & 3     & 1750  & 500   & 52    & 0.1\\
        22  & 13    & 7     & 7     & 2     & 1550  & 800   & 32    & 0.07\\
        23  & 12    & 8     & 10    & 8     & 1800  & 750   & 38    & 0.16\\
        24  & 11    & 13    & 5     & 1     & 1500  & 800   & 44    & 0.11\\
        25  & 6     & 15    & 3     & 2     & 1700  & 550   & 66    & 0.12\\
        \hline\hline
    \end{tabular}
    \label{tab: validation study parameters 1}
\end{table}

\begin{table}[hbt!]
    \centering
    \caption{Parameters for validation study -- 2}
    \begin{tabular}{lrrrrrrrrrr}
        \hline\hline
        Inst. \# & $N_{\text{plane}}$ & $N_{\text{parking}}$ & $N_{\text{sats}}$ & $R_1$ & $R_2$ & $Q_1$ & $Q_2$ & $k_{R_1,3}$ & $k_{Q_1,3}$ & $\alpha_{w}$ \\\hline
        1   & 22    & 7     & 60    & 3     & 2     & 7     & 10    & 10    & 17    & 1.4\\
        2   & 36    & 11    & 40    & 3     & -1    & 18    & 6     & 8     & 13    & 0.4\\
        3   & 22    & 3     & 60    & 8     & 7     & 8     & 5     & 8     & 16    & 0.8\\
        4   & 30    & 20    & 35    & 2     & -2    & 2     & 3     & 7     & 18    & 1.8\\
        5   & 28    & 16    & 40    & 5     & -2    & 10    & 7     & 2     & 6     & 0.9\\
        6   & 32    & 5     & 20    & 1     & 0     & 3     & 8     & 5     & 13    & 0.7\\
        7   & 20    & 17    & 45    & 4     & 3     & 4     & 6     & 1     & 3     & 1.7\\
        8   & 40    & 12    & 20    & 2     & -1    & 14    & 6     & 11    & 20    & 1\\
        9   & 20    & 12    & 45    & 17    & 0     & 7     & 7     & 14    & 10    & 1.9\\
        10  & 32    & 18    & 55    & 11    & 2     & 3     & 5     & 4     & 19    & 0.2\\
        11  & 34    & 19    & 45    & 15    & 6     & 14    & 2     & 3     & 15    & 1.2\\
        12  & 34    & 7     & 35    & 18    & 9     & 20    & 2     & 17    & 10    & 1.5\\
        13  & 30    & 17    & 30    & 15    & 7     & 19    & 2     & 7     & 1     & 1.5\\
        14  & 32    & 11    & 25    & 10    & 8     & 9     & 8     & 15    & 3     & 1.6\\
        15  & 36    & 4     & 50    & 6     & 3     & 12    & 3     & 18    & 14    & 0.3\\
        16  & 30    & 20    & 25    & 13    & 6     & 11    & 4     & 9     & 1     & 2\\
        17  & 24    & 9     & 50    & 12    & 5     & 15    & 9     & 12    & 2     & 0.5\\
        18  & 38    & 10    & 30    & 9     & 4     & 5     & 1     & 16    & 8     & 1.2\\
        19  & 38    & 2     & 55    & 20    & 10    & 17    & 10    & 13    & 12    & 1.3\\
        20  & 24    & 6     & 25    & 16    & 5     & 1     & 4     & 20    & 7     & 0\\
        21  & 26    & 8     & 40    & 19    & 1     & 13    & 1     & 4     & 11    & 0.1\\
        22  & 26    & 1     & 50    & 13    & 10    & 18    & 3     & 6     & 4     & 0.6\\
        23  & 22    & 14    & 30    & 12    & 9     & 10    & 9     & 13    & 9     & 1.4\\
        24  & 28    & 13    & 35    & 7     & 1     & 16    & 10    & 19    & 2     & 0.6\\
        25  & 40    & 15    & 22    & 14    & 8     & 6     & 8     & 18    & 5     & 1.1\\
        \hline\hline
    \end{tabular}
    \label{tab: validation study parameters 2}
\end{table}

\end{document}